\documentclass[11pt]{amsart}
\usepackage{amssymb}
\usepackage{verbatim}

\newtheorem{lemma}{Lemma} [section]
\newtheorem{thm}[lemma]{Theorem}

\newtheorem{cor}[lemma]{Corollary}
\newtheorem{prop}[lemma]{Proposition}

\theoremstyle{remark}

\newtheorem*{remark}{Remark}

\theoremstyle{definition}

\numberwithin{equation}{section}

\DeclareMathOperator{\lcm}{{lcm}}
\newcommand{\C}{{\mathbb C}}
\newcommand{\bC}{\C}
\newcommand{\bZ}{{\mathbb Z}}

\newcommand{\cL}{{\mathcal L}}
\newcommand{\tth}{^{\operatorname{th}}}

\newcommand{\Line}{{\mathbb P}^1}
\renewcommand{\hat}{\widehat}

\begin{document}



\title{Decompositions of Laurent polynomials}

\author{Michael E. Zieve}
\address{
  Center for Communications Research,
  805 Bunn Drive,
  Princeton, NJ 08540
}
\email{zieve@math.rutgers.edu}
\urladdr{www.math.rutgers.edu/$\sim$zieve/}

\thanks{We thank Fedor Pakovich for informing us of his papers applying
Ritt's results, and Yuri Bilu for suggesting (via Pakovich) that
decompositions of Laurent polynomials are
related to genus-zero factors of curves $g(x)-h(y)$ where $g,h\in\C[x]$
have $\gcd(\deg(g),\deg(h))=2$.}

\begin{abstract}
In the 1920's, Ritt studied the operation of functional composition
$g\circ h(x) = g(h(x))$ on complex rational functions.  In the case
of polynomials, he described all the ways in which a polynomial can
have multiple `prime factorizations' with respect to this operation.
Despite significant effort by Ritt and others, little progress has
been made towards solving the analogous problem for rational functions.
In this paper we use results of Avanzi--Zannier and Bilu--Tichy to prove
analogues of Ritt's results for decompositions of Laurent polynomials,
i.e., rational functions with denominator $x^n$.
\end{abstract}


\maketitle


\section{Introduction}
In the 1920's, Ritt \cite{Ritt} studied the possible ways of writing a complex
polynomial as a composition of lower-degree polynomials.  To this end,
a polynomial $f\in\C[x]$ with $\deg(f)>1$ is called \emph{indecomposable}
if it cannot be written as a composition $f(x)=g(h(x))$ with $g,h\in\C[x]$
and $\deg(g),\deg(h)<\deg(f)$.  By induction, any polynomial of degree more
than one
can be written as the composition of indecomposable polynomials.  Although
this decomposition need not be unique, Ritt proved that its length is
unique, and moreover he gave a recursive procedure for obtaining any
decomposition from any other.  Ritt's results are quite fundamental, and
have been applied in various wide-ranging contexts (cf.\ \cite{BWZ,BT,GTZ,GTZ2,
MP,P1,P2,PRY,Zannier2}, among others).

Unfortunately, there are no known analogues of Ritt's results in the case
of rational functions.  Ritt himself was the first to study this
\cite{Rittrat,Ritt2}.  He noted \cite{Ritt2} that the action of the group
$A_4$ on
the Riemann sphere, together with the fact that $A_4$ has maximal chains
of subgroups $1<C_2<V_4<A_4$ and $1<C_3<A_4$, implies that a certain
degree-$12$ rational function can be written as both the composition
of two indecomposables and the composition of three indecomposables.
(This example is reproduced in the context of modular forms in \cite{GS,MS}.)
Further, if $f(x)$ is the map on $x$-coordinates induced by
multiplication-by-$p$ on the elliptic curve $y^2=x^3+1$, for any prime
$p$ with $p\equiv 2\pmod{3}$, then $f$ is indecomposable but there is
a decomposable $g\in\C(x)$ for which $x^3\circ f = g\circ x^3$ \cite{LZ}.
Further families of counterexamples to the rational function
analogues of Ritt's results are given in \cite{LZ}; however, as noted there,
all known examples fit into one of three
simple types, which suggests there may be a concise description of all
examples.  On the other hand, proving such a possibility seems far beyond
current techniques.

In this paper we study a situation which lies between the polynomial and
rational function cases: namely, we study \emph{Laurent polynomials}, i.e.,
rational functions of the form $f(x)/x^n$ with $f\in\C[x]$.  We will prove
that decompositions of Laurent polynomials satisfy variants of Ritt's results.
Our statements involve the Dickson polynomials $D_n(x)$, which are defined by
the functional equation $D_n(x+1/x)=x^n+1/x^n$; these are related to the
classical Chebychev polynomials $T_n(x)$ via $D_n(x)=2T_n(x/2)$.
We say a rational function of degree $>1$ is
\emph{indecomposable} if it cannot be written as the composition of
rational functions of strictly lower degrees, and a
\emph{complete decomposition} of a rational function is an expression
of the rational function as the composition of indecomposable rational
functions.  We note (cf.\ Lemma~\ref{fields}) that a decomposable Laurent
polynomial can actually be written as the composition of two Laurent
polynomials of strictly lower degrees, rather than just as the composition
of lower-degree rational functions.
Writing $\cL$ for the set of all complex Laurent polynomials,
our Laurent polynomial analogue of the classical `first theorem of Ritt' is
as follows:

\begin{thm}
\label{LRitt1}
If $f=p_1\circ p_2\circ\dots\circ p_r=q_1\circ q_2\circ\dots\circ q_s$ where
$p_i,q_j\in\C(x)$ are indecomposable and $f\in\cL$, then
the sequences $(\deg(p_1),\dots,\deg(p_r))$ and
$(\deg(q_1),\dots,\deg(q_s))$ are permutations of one another (so $r=s$).
Moreover, there is a finite sequence
of complete decompositions of $f$ which begins with
$p_1\circ\dots\circ p_r$ and ends
with $q_1\circ\dots\circ q_s$, where consecutive decompositions in the
sequence differ only in that two adjacent indecomposables in the first
decomposition are replaced in the second decomposition by two others having
the same composition.
\end{thm}

Our Laurent polynomial analogue of the `second theorem of Ritt' is:
\begin{thm}
\label{LRitt2}
If $f=g_1\circ h_1 = g_2\circ h_2$ where $g_1,g_2,h_1,h_2\in\C(x)$ are
indecomposable and $f\in\cL$, then (after perhaps exchanging the pairs
$(g_1,h_1)$ and $(g_2,h_2))$ there exist degree-one
$\mu_1,\dots,\mu_4\in\C(x)$ such that
\begin{align*}
g_1 &= \mu_1\circ G_1\circ\mu_3 \\
g_2 &= \mu_1\circ G_2\circ\mu_4 \\
h_1 &= \mu_3^{-1}\circ H_1\circ\mu_2 \\
h_2 &= \mu_4^{-1}\circ H_2\circ\mu_2,
\end{align*}
where one of the following holds (with $n$ prime):
\begin{enumerate}
\item[(\thethm.1)] $G_1=G_2$ and $H_1=H_2$ with $G_1,H_2\in\cL$ and
either $G_1\in\C[x]$ or $H_2=x^n$;
\item[(\thethm.2)] $G_1=H_2=x^n$, $\,H_1=x^r q(x^n)$, and
$G_2=x^r q(x)^n$ with $q\in\C(x)$ and $r\in\bZ_{>0}$ coprime to $n$;
\item[(\thethm.3)] $G_1=H_2=D_m$ and $H_1=G_2=D_n$,
where $m\ne n$ is prime;
\item[(\thethm.4)] $G_1=D_n$, $\,H_1=G_2=x+1/x$, and
$H_2=x^n$;
\item[(\thethm.5)] $G_1=G_2=D_n$, $\,H_1=x+1/x$, and
$H_2=\zeta x + 1/(\zeta x)$, where $\zeta^n=1$.
\end{enumerate}
\end{thm}

We emphasize that, in (\ref{LRitt2}.2), we do not require
$q\in\cL$.  In fact, our proof shows we can require either
$q\in\cL$ or $q=Q(\frac{1}{x+1})$ with $Q\in x\C[x]$.  To see why
the latter case gives rise to Laurent polynomials (after composing
with $\mu_2$),
put $q=Q(\frac{1}{x+1})$ with $Q\in x\C[x]$, so
$xq(x^2)\circ i\frac{x-1}{x+1} = i\frac{x-1}{x+1}Q(\frac{(x+1)^2}{4x})$, which
is in $\cL$.

These results generalize the classical theorems of Ritt, which are
obtained by requiring all the rational functions to be polynomials.
Stated in the other direction, if we begin with Ritt's results and
attempt to generalize them to decompositions of Laurent polynomials, we
must replace the various polynomials in Ritt's results by
rational functions, and also we must allow the new possibilities
(\ref{LRitt2}.4) and (\ref{LRitt2}.5).  In fact, (\ref{LRitt2}.5)
can be obtained from two applications of (\ref{LRitt2}.4), in addition
to composing with linears: for, if $\zeta^n=1$ then
\begin{align*}
D_n \circ \left(\zeta x + \frac{1}{\zeta x}\right) &=
 D_n\circ \left(x+\frac{1}{x}\right)\circ
  \zeta x \\
&= \left(x+\frac{1}{x}\right)\circ x^n \circ \zeta x \\
&= \left(x+\frac{1}{x}\right)\circ x^n \\
&= D_n\circ \left(x + \frac{1}{x}\right).
\end{align*}

One consequence of Ritt's results, which actually was deduced as a step
in Ritt's proofs, is a certain `rigidity' property of polynomial
decompositions:

\begin{cor}
\label{unique}
If $g_1\circ h_1=g_2\circ h_2$ where $g_1,g_2,h_1,h_2\in\C[x]\setminus\C$ and\/
$\deg(g_1)=\deg(g_2)$, then there is
a linear $\mu\in\C[x]$ such that $g_2=g_1\circ\mu$ and $h_2=\mu^{-1}\circ h_1$.
\end{cor}

Note that (\ref{LRitt2}.5) provides counterexamples to the Laurent polynomial
analogue of Corollary~\ref{unique}.  Further counterexamples are obtained
by putting $n=2$ in (\ref{LRitt2}.4).  We will determine all examples:

\begin{prop}
\label{Lunique}
If $f=g_1\circ h_1=g_2\circ h_2$ where $f\in\cL\setminus\C$ and
$g_1,g_2,h_1,h_2\in\C(x)$
satisfy\/ $\deg(g_1)=\deg(g_2)$, then, perhaps after exchanging
$(g_1,h_1)$ and $(g_2,h_2)$, there exist $G\in\C[x]$, $H\in\cL$, and
degree-one $\mu_1,\mu_2\in\C(x)$ such that
\begin{align*}
g_1&=G\circ G_1\circ\mu_1 \\
g_2&=G\circ G_2\circ\mu_2 \\
h_1&=\mu_1^{-1}\circ H_1\circ H \\
h_2&=\mu_2^{-1}\circ H_2\circ H,
\end{align*}
where one of the following holds (in which $n\in\bZ_{>0}$):
\begin{enumerate}
\item[(\thethm.1)]$G_1=G_2=H_1=H_2=x$;
\item[(\thethm.2)]$G_1=H_2=x^n$, $\,H_1=(x^n+1)/x^r$, and
$G_2=(x+1)^n/x^r$, where $0< r < n$ and $\gcd(r,n)=1$;
\item[(\thethm.3)]$G_1=-G_2=D_n$, $\,H_1=x+1/x$, and
 $H_2=\zeta x+1/(\zeta x)$, where $\zeta^n=-1$;
\item[(\thethm.4)]$G_1=D_2$, $\,H_1=G_2=x+1/x$, and
 $H_2=x^2$.
\end{enumerate}
Moreover, in \emph{(\thethm.2)}--\emph{(\thethm.4)} we may assume
$H=\alpha x^s$ with $\alpha\in\C^*$ and $s\in\bZ_{>0}$.
\end{prop}

Ritt proved a generalization of the polynomial version of
Theorem~\ref{LRitt2}, which can
be used to describe all polynomials $g_1,g_2,h_1,h_2$ with
$g_1\circ h_1 = g_2\circ h_2$ \cite{BWZ}.  We will prove the following
analogue for Laurent polynomials:

\begin{thm} \label{Lbidec}
Let $f\in\cL\setminus\C$ and $g_1,g_2,h_1,h_2\in\bC(x)$ satisfy
$f=g_1\circ h_1 = g_2\circ h_2$.  Then, perhaps after switching
$(g_1,h_1)$ and $(g_2,h_2)$, we have
\begin{align*}
g_1&=G\circ G_1\circ\mu_1\\
g_2&=G\circ G_2\circ\mu_2\\
h_1&=\mu_1^{-1}\circ H_1\circ H\\
h_2&=\mu_2^{-1}\circ H_2\circ H
\end{align*}
for some $G\in\C[x]$, some $H\in\cL$, and some degree-one
$\mu_1,\mu_2\in\C(x)$, where one of the following holds
(in which $m,n$ are coprime positive integers, and $p\in\C[x]\setminus\{0\}$):
\begin{enumerate}
\item[(\thethm.1)]$G_1=H_2=x^n$, $\,H_1=x^r p(x^n)$,
and $G_2=x^r p(x)^n$, where $r\in\bZ$ with $\gcd(r,n)=1$;
\item[(\thethm.2)]$G_1=x^2$, $\,H_1=(x-\frac{1}{x})p(x+\frac{1}{x})$,
$\,G_2=(x^2-4)p(x)^2$, and $H_2=x+1/x$;
\item[(\thethm.3)]$G_1=H_2=D_m$ and $H_1=G_2=D_n$;
\item[(\thethm.4)]$G_1=(\frac{x^2}{3}-1)^3$, $\,H_1=x^2+2x+\frac{1}{x}-
\frac{1}{4x^2}$,
$\,G_2=3x^4-4x^3$, and $H_2=\frac{1}{3}((x+1-\frac{1}{2x})^3+4)$;
\item[(\thethm.5)]$G_1=D_{dm}$, $\,H_1=x^n+1/x^n$,
$\,G_2=-D_{dn}$, and $H_2=(\zeta x)^m+1/(\zeta x)^m$, where
$d\in\bZ_{>1}$ and $\zeta^{dmn}=-1$;
\item[(\thethm.6)] $G_1=D_m$, $\,H_1=G_2=x^n+1/x^n$,
and $H_2=x^m$.
\end{enumerate}
Moreover, in all cases besides \emph{(\thethm.1)} and \emph{(\thethm.3)},
we may assume $H=\alpha x^s$ with $\alpha\in\C^*$ and $s\in\bZ_{>0}$.
\end{thm}

The analogous result for decompositions of polynomials \cite{BWZ}
involves only cases (\ref{Lbidec}.1) and (\ref{Lbidec}.3).

Ritt's proofs of the polynomial versions of Theorems \ref{LRitt1} and
\ref{LRitt2} are independent of one another, and have quite distinct
flavors.  His proof of Theorem \ref{LRitt1} for polynomials is essentially
group theoretic: if $f$ is a polynomial then the inertia group $I$ at any
infinite place of (the Galois closure of) $\C(x)/\C(f(x))$ is transitive,
so one can translate questions about decompositions of $f$ into
questions about subgroups of $I$, which are not difficult to resolve since
$I$ is cyclic.  On the other hand, Ritt's proof of Theorem \ref{LRitt2}
for polynomials is a genus computation, as he determines all polynomials
$g_1,h_1$ of coprime degrees for which the curve $g_1(x)-h_1(y)$ has
genus zero.  For Laurent polynomials we require a different approach,
since there is no longer a transitive inertia group, so Theorem \ref{LRitt1}
cannot be proved via group theory.  Instead we first prove
Theorem~\ref{Lbidec}, using results of Avanzi--Zannier \cite{AZ} and
Bilu--Tichy \cite{BT}, which in turn rely on Ritt's second theorem and
related genus computations (among other things).
After determining the possible decompositions of the specific rational
functions appearing in Theorem~\ref{Lbidec}, we can then deduce
Theorems~\ref{LRitt1} and \ref{LRitt2}.  We pay special attention to
decompositions of $H_1$ and $G_2$ from (\ref{LRitt2}.2), in view of their
role in potential analogues of Ritt's results for rational functions:
these $H_1$ and $G_2$ are especially important since they have the same
shape as one of the main sources of rational function counterexamples
(the one including the elliptic curve examples mentioned above).

Ritt's proofs used the language of Riemann surfaces; several authors have
rewritten his proofs in different languages
\cite{Binder,DW,Engstrom,Fried-Ritt,FriedMacRae,LN,Levi,McConnell,Mueller,
Schinzelold,Schinzel,Tortrat,Zannier}.  For some applications the recursive
procedure in Theorem~\ref{LRitt2} is not sufficient, and one needs more
precise information about the collection of all the different decompositions
of a polynomial; see \cite{MZ} for the state of the art on polynomial
decomposition.  We do not know whether there are Laurent polynomial
analogues of the latter results.

The contents of this paper are as follows.  In the next section we
prove some general results about decompositions of Laurent polynomials,
based on which we outline our strategy for proving our main results.
In Sections \ref{sec-easy} and \ref{sec-hard} we describe all decompositions
of the various special Laurent polynomials occurring in the statements
of the above results.  We use these specific decompositions to
prove preliminary versions of Theorem~\ref{Lbidec} in Sections \ref{sec-1} and
\ref{sec-2}, and finally we conclude in Section \ref{sec-proofs} by
proving the results stated in this introduction.


\section{Preliminary reductions}

Recall that the set $\cL$ of Laurent polynomials consists of all
rational functions whose denominator is a power of $x$, or equivalently,
all rational functions having no poles besides $0$ and $\infty$.  This
perspective yields the following result:

\begin{lemma}\label{basic}
If $f=g\circ h$ where $f\in\cL\setminus\C$ and $g,h\in\C(x)$,
then there is a degree-one $\mu\in\C(x)$
such that $G:=g\circ\mu$ and $H:=\mu^{-1}\circ h$ satisfy
one of the following:
\begin{enumerate}
\item[(\thethm.1)] $G\in\C[x]$ and $H\in\cL$;
\item[(\thethm.2)] $G\in\cL$ and $H = x^n$
for some $n\in\bZ_{>0}$.
\end{enumerate}
\end{lemma}

\begin{proof}
The poles of $f=g\circ h$ are the preimages under $h$ of the poles
of $g$; by hypothesis, these preimages form a subset of $\{0,\infty\}$.
Hence $g$ has at most two poles.  First suppose $g$ has a unique pole,
say $\alpha$.  Pick a degree-one $\mu\in\C(x)$ for which
$\mu(\infty)=\alpha$, so
that $G:=g\circ\mu$ has $\infty$ as its unique pole, whence
$G\in\C[x]$.  Then $f=G\circ H$ where
$H:=\mu^{-1}\circ h$, and $H$ can have no poles besides
$0$ and $\infty$, so $H\in\cL$, as in (\ref{basic}.1).
Now suppose $g$ has two poles, say $\alpha$ and $\beta$.  Since
$g\circ h$ has at most two poles, both $\alpha$ and $\beta$ must have
unique preimages under $h$, which must be $0$ and $\infty$.  Say
$\alpha=h(0)$ and $\beta=h(\infty)$, and put $\gamma=h(1)$.
Pick a degree-one $\mu\in\C(x)$ which maps $0\mapsto\alpha$ and
$\infty\mapsto\beta$ and $1\mapsto\gamma$.
Then the poles of $G:=g\circ\mu$ are $0$ and $\infty$, so $G\in\cL$,
and $H:=\mu^{-1}\circ h$ has its unique pole at $\infty$
(so $H\in\C[x]$) and has $0$ as its unique root (so $H$ is a
monomial) and maps $1\mapsto 1$ (so $H$ is monic).
\end{proof}

Thus, in what follows we will restrict to decompositions
$f=G\circ H$ where $G$ and $H$ satisfy (\ref{basic}.1)
or (\ref{basic}.2).
We refer to decompositions of these types as `Type 1' and `Type 2'
decompositions.  A pair of decompositions of the same Laurent polynomial
must be in one of three categories: both decompositions could be Type 1,
both could be Type 2, or one could be Type 1 and the other Type 2.
It is easy to describe the pairs of Type~2 decompositions of a Laurent
polynomial:

\begin{prop}\label{Type2}
If $g_1\circ x^n = g_2\circ x^m$ with $g_i\in\cL$ and $n,m>0$,
then there exists $G\in\cL$ such that
$g_1=G\circ x^{\lcm(n,m)/n}$ and $g_2=G\circ x^{\lcm(n,m)/m}$.
\end{prop}

In other words, if we write a Laurent polynomial $f$ as $f=G\circ x^N$
with $N$ maximal, then every Type 2 decomposition of $f$ is (up to linears)
$G(x^n)\circ x^{N/n}$.

\begin{proof}
Writing $f=g_1\circ x^n$, the field $\C(f)$ is contained in
$\C(x^n)\cap\C(x^m)=\C(x^d)$, where $d=\lcm(n,m)$.  Write $d=Nn=Mm$,
so $g_1\circ x^n = G_1\circ x^d$ for some $G_1\in\C(x)$ (which is
automatically a Laurent polynomial), whence $g_1 = G_1\circ x^N$.
Likewise $g_2 = G_2\circ x^M$, and we have $G_1\circ x^d = f =
G_2\circ x^d$, so $G_1=G_2$.  Thus $f=G_1(x^{Nn})$, and its two Type 2
decompositions are $G_1(x^N)\circ x^n$ and $G_1(x^M)\circ x^m$.
\end{proof}

Next we consider Laurent polynomials with two Type 1 decompositions:
$f=g_1\circ h_1=g_2\circ h_2$ with $g_i\in\C[x]$ and $h_i\in\cL$.
Then there is an irreducible factor
$E(x,y)$ of $g_1(x)-g_2(y)$ such that $E(h_1(x),h_2(x))=0$, so $E(x,y)=0$
defines a genus-zero curve having at most two closed points lying over
$x=\infty$ (since $f$ has at most two poles).
To classify the possibilities in this case, we
use a result of Bilu and Tichy \cite{BT} describing the polynomials
$g_1,g_2$ for which the curve $g_1(x)=g_2(y)$ has an irreducible component
with these properties.  Note that in this situation there automatically
exist nonconstant $h_1,h_2\in\cL$ such that $g_1\circ h_1=g_2\circ h_2$,
coming from a rational parametrization of the component in question.

Finally we consider Laurent polynomials with decompositions of both types:
$f=g_1\circ h_1 = g_2\circ x^n$ where $g_1\in\C[x]$ and $h_1,g_2\in\cL$
(and $n>1$).  Letting $\zeta$ be a primitive $n\tth$ root of unity,
we have
\[
g_1\circ h_1(\zeta x) = g_2\circ x^n\circ \zeta x = g_2\circ x^n=
g_1\circ h_1(x).\]
Let $h_2(x)=h_1(\zeta x)$.
To classify the possibilities where $h_2\ne h_1$, we use a result
of Avanzi and Zannier \cite{AZ} describing the polynomials $g_1$ for which
there are distinct nonconstant rational functions $h_1,h_2$ such that
$g_1\circ h_1=g_1\circ h_2$.  Finally, if $h_1(\zeta x)=h_1(x)$ then
$h_1=H(x^n)$
for some $H\in\cL$, where $g_1\circ H=g_2$.  Thus, these possibilities come
from decompositions of the Laurent polynomial $g_2$, which can be
controlled inductively.

We now recall the well-known connection between decompositions of a
rational function $f$ and intermediate fields between $\C(x)$ and $\C(f(x))$,
as well as the corresponding results for polynomials and Laurent polynomials.

\begin{lemma}\label{fields}
For $f\in\C(x)\setminus\C$, the fields between\/ $\C(x)$ and\/ $\C(f)$ are
precisely the fields\/ $\C(h)$, where $g,h\in\C(x)$ satisfy $f=g\circ h$;
moreover, for $h,H\in\C(x)$, we have\/ $\C(h)=\C(H)$ if and only if there
is a degree-one $\mu\in\C(x)$ such
that $h=\mu\circ H$.  If $f$ is a Laurent polynomial (respectively,
polynomial) and $f=g\circ h$
with $g,h\in\C(x)$, then there is a degree-one $\mu\in\C(x)$ such that
both $g\circ\mu$ and $\mu^{-1}\circ h$ are Laurent polynomials (respectively,
polynomials).
\end{lemma}

\begin{proof}
The first statement follows from L\"uroth's theorem.  Now suppose $f=g\circ h$
where $g,h\in\C(x)$ and $f\in\C[x]$; since $\infty$ is the unique pole of $f$,
it follows that $g$ has a unique pole $\alpha$, and $\infty$ is the unique
preimage of $\alpha$ under $h$.  Pick a degree-one $\mu\in\C(x)$ which maps
$\infty\mapsto\alpha$, so both $g\circ\mu$ and $\mu^{-1}\circ h$ are
rational functions whose unique pole is $\infty$, hence they are polynomials.
Next suppose $f=g\circ h$ where $g,h\in\C(x)$ and $f\in\cL$; then $f$ has no
poles besides $0$ and $\infty$, so $g$ also has at most two poles, and the
preimages of these poles under $h$ are a subset of $\{0,\infty\}$.
Pick a degree-one $\mu\in\C(x)$ which maps the poles of $g$ to either
$\{\infty\}$ or $\{0,\infty\}$; then both $g\circ\mu$ and $\mu^{-1}\circ h$
have no poles outside $\{0,\infty\}$, hence are Laurent polynomials.
\end{proof}


\section{Decompositions of Laurent polynomials of special types}
\label{sec-easy}

In this section we describe all decompositions of certain special
Laurent polynomials occurring in our results.  Knowledge of these
decompositions will be used in the proofs of our main results.

We begin with $f=x^n+1/x^n$ (where $n\in\bZ_{>0}$), whose decompositions
turn out to be the main source of Laurent polynomial decompositions
that are not polynomial decompositions.

\begin{lemma}\label{Dickson}
If $g,h\in\C(x)$ satisfy $g\circ h=x^n+x^{-n}$ for some $n>0$,
then there is a divisor $d$ of $n$ and a degree-one $\mu\in\C(x)$ such that
one of the following holds:
\begin{enumerate}
\item[(\thethm.1)]$g\circ\mu = x^{n/d}+x^{-n/d}$ and
$\mu^{-1}\circ h = x^d$;
\item[(\thethm.2)]$g\circ\mu = \beta^n D_{n/d}$ and
$\mu^{-1}\circ h = (x/\beta)^d+(\beta/x)^d$ where $\beta^{2n}=1$.
\end{enumerate}
\end{lemma}

\begin{proof}
Writing $f=x^n+x^{-n}$, we see that $\C(x)/\C(f)$ is
Galois, with Galois group $G$ being dihedral of order $2n$ and consisting
of the automorphisms $x\mapsto \zeta x^e$ with $\zeta^n=1$ and $e\in\{1,-1\}$.
Let $C$ be the cyclic subgroup of $G$ consisting of the automorphisms
$x\mapsto\zeta x$.  Let $H$ be a subgroup of $G$, and let $d=\#(H\cap C)$;
then $H\cap C$ consists of the automorphisms $x\mapsto\delta x$ with
$\delta^d=1$, so the fixed field $\C(x)^{H\cap C}$ equals $\C(x^d)$.
If $H=H\cap C$ then the chain of groups $1<H<G$ corresponds (via 
Lemma~\ref{fields}) to the
decomposition $f=(x^{n/d}+x^{-n/d})\circ x^d$.  Now suppose $H\ne H\cap C$,
so $\#H=2d$.  Pick some $\zeta$ for which $H$ contains the automorphism
$x\mapsto \zeta/x$.  Then $\C(x)^H=\C(x^d+(\zeta/x)^d)=
\C((x/\beta)^d+(\beta/x)^d)$
where $\beta^2=\zeta$ (so $\beta^{2n}=1$), and the corresponding decomposition
is $f=(\beta^n D_{n/d})\circ ((x/\beta)^d+(\beta/x)^d)$.
\end{proof}

We also recall the possible decompositions of $x^n$ and $D_n$:
\begin{lemma}
If $g\circ h=x^n$ with $g,h\in\C[x]$ and $n>0$, then there is a linear
$\mu\in\C[x]$ and a divisor $d$ of $n$ such that $g\circ\mu=x^d$ and
$\mu^{-1}\circ h=x^{n/d}$.  If $g\circ h=D_n$ with $g,h\in\C[x]$ and
$n>0$, then there is a linear $\mu\in\C[x]$ and a divisor $d$ of $n$
such that $g\circ\mu=D_d$ and $\mu^{-1}\circ h=D_{n/d}$.
\end{lemma}

\begin{proof}
This follows from Corollary~\ref{unique}, together with the fact that
$D_d\circ D_{n/d}=D_n$ (which follows from the functional equation defining
$D_n$).
\end{proof}

Rather than writing out all the decompositions of the rational functions in
(\ref{Lbidec}.4), we show that (\ref{Lbidec}.4) is a consequence of
(\ref{Lbidec}.1) and (\ref{Lbidec}.2), if we allow compositions with
linear polynomials.  Namely, putting
$p=\frac{x}{2}+\sqrt{2}$ and $\nu=x\sqrt{2}$, we have
\[
x^2+2x+\frac{1}{x}-\frac{1}{4x^2} =
 \left(x+\frac{1}{x}\right)\cdot p\left(x-\frac{1}{x}\right)\circ\nu,\]
so for
\[
f:=\left(\frac{x^2}{3}-1\right)^3\circ
 \left(x^2+2x+\frac{1}{x}-\frac{1}{4x^2}\right)\]
we have
\begin{align*}
f &= \left(\frac{x}{3}-1\right)^3 \circ x^2\circ
        \left(x+\frac{1}{x}\right)\cdot p\left(x-\frac{1}{x}\right)\circ\nu \\
&= \left(\frac{x}{3}-1\right)^3 \circ (x^2+4)p(x)^2\circ
      \left(x-\frac{1}{x}\right)\circ\nu,
\end{align*}
where the last equality comes from (\ref{Lbidec}.2).  Now put
$\mu=\sqrt{2}(x-1)$, so
\[
(x^2+4)p(x)^2\circ\mu = x^4+4x+3\] and
\[\mu^{-1}\circ \left(x-\frac{1}{x}\right)\circ\nu = x+1-\frac{1}{2x},\]
and thus if we put $\lambda=3x-4$ then
\begin{align*}
f &= x^3 \circ \left(\frac{x}{3}-1\right)\circ (x^4+4x+3)\circ
 \left(x+1-\frac{1}{2x}\right) \\
&= x^3 \circ \frac{x^4+4x}{3} \circ \left(x+1-\frac{1}{2x}\right) \\
&= x\left(\frac{x+4}{3}\right)^3 \circ \lambda\circ\lambda^{-1}\circ x^3 \circ
       \left(x+1-\frac{1}{2x}\right) \quad\text{(from (\ref{Lbidec}.1))}\\
&= (3x-4)x^3\circ \frac{x+4}{3}\circ x^3 \circ \left(x+1-\frac{1}{2x}\right) \\
&= (3x^4-4x^3)\circ \frac{(x+1-\frac{1}{2x})^3+4}{3}.
\end{align*}


\section{Decompositions of Ritt-twistable Laurent polynomials}
\label{sec-hard}

In this section we study decompositions of the Laurent polynomials occurring
in (\ref{Lbidec}.1) and (\ref{Lbidec}.2).  Some of the results we prove
will be used in the proofs of our main results.  We also prove other
results giving a full picture of the decompositions of these special
Laurent polynomials, in view of the important role these examples play
in the study of rational function analogues of Ritt's results.

Case (\ref{Lbidec}.1) involves
Laurent polynomials of the form $x^r q(x^n)$ and $x^r q(x)^n$, where
$q\in\cL\setminus\{0\}$ and $\gcd(r,n)=1$.  These are the natural Laurent
polynomial analogues of the polynomials occurring in Ritt's results
(which have the same shape but with $q\in\C[x]$).
The Laurent polynomials in (\ref{Lbidec}.2), however, have a different
shape, namely
$H_1=(x-1/x)p(x+1/x)$ and $G_2=(x^2-4)p(x)^2$, with
$p\in\C[x]\setminus\{0\}$.  We now show that there are linear changes of
variables which transform $H_2$ and $G_2$ into the same general
shape as the previous Laurent polynomials, namely $x q(x^2)$ and $x q(x)^2$,
although here we must allow $q$ to be a rational function that is not in
$\cL$.  Specifically, if we put
\stepcounter{lemma}
\begin{equation}\label{q}\tag{\thethm.1}
q=4i\frac{p(2\frac{x-1}{x+1})}{x+1},\end{equation}
then
\begin{align}
\label{q2}\tag{\thethm.2}x q(x^2) &= H_1 \circ \frac{x+i}{x-i} \\
\label{q3}\tag{\thethm.3}x q(x)^2 &= G_2 \circ \frac{2x-2}{x+1}.
\end{align}
It is shown in \cite{MZ} that a polynomial of the form $x^r q(x^n)$
(with $\gcd(r,n)=1$) can only decompose into polynomials of
the same shape (composed with linears), and likewise for $x^r q(x)^n$.
We will prove the analogous result for Laurent polynomials;
the corresponding assertion is not generally true when $q$ is one of
the rational functions in (\ref{q}), but nevertheless we determine
all decompositions in this situation.
We remark (cf.\ \cite{LZ}) that Ritt's original $A_4$ example (after
linear changes) provides an example of an `odd' rational function
$x q(x^2)$ which can be written as the composition of two rational
functions that are not linear changes of odd rational functions;
similar examples occur for $q$ as in (\ref{q}).

\begin{prop}\label{Ritttwisteasy}
Let $n,r\in\bZ$ satisfy $n>1$ and $\gcd(n,r)=1$, and pick $p\in\C[x]$
with $x\nmid p$.  Suppose $g,h\in\C(x)$ satisfy $g\circ h = x^r p(x)^n$.
Then there is a degree-one $\mu\in\C(x)$ such that
$g\circ\mu=x^i G^n$ and
$\mu^{-1}\circ h=x^j H^n$ for some $\delta\in\C$, some $i,j\in\bZ$,
and some $G,H\in\C[x]$.
\end{prop}

\begin{proof}
If $r\ge 0$ then $x^r p(x)^n$ is a polynomial, in which case the result
is proved in \cite{MZ} if $g,h\in\C[x]$, and the general case follows
from Lemma~\ref{fields}.  Henceforth assume $r<0$.

By Lemma~\ref{basic}, after replacing $g$ and $h$ by $g\circ\mu$ and
$\mu^{-1}\circ h$ for suitable degree-one $\mu\in\C(x)$, we may assume
$g,h\in\cL$ and either $g\in\C[x]$ or $h=x^m$ with $m\in\bZ_{>0}$.
First suppose $h=x^m$.  Letting $\zeta$ be a primitive $m\tth$ root of
unity, we have $g\circ h(\zeta x)=g\circ h(x)$, so
$\zeta^r x^r p(\zeta x)^n = x^r p(x)^n$.  Thus $p(\zeta x)$ is a constant
times $p(x)$, so $p=x^s G(x^m)$ with $G\in\C[x]$ and $x\in\bZ_{\ge 0}$.
Since $x^r p(x)^n=g\circ x^m$, we have $r+ns=mi$ with
$i\in\bZ_{\ge 0}$, so $g=x^i G(x)^n$.  Putting $j=m$ and $H=1$ gives
the desired conclusion.  Henceforth assume $g\in\C[x]$.

Write $h = A/x^s$ where $s\in\bZ_{>0}$ and $A\in\C[x]$ with $x\nmid A$.
Write $g=\theta\prod_{\alpha} (x-\alpha)^{n_\alpha}$, where the $\alpha$
are the distinct complex roots of $g$ (and $n_\alpha\in\bZ_{>0}$ and
$\theta\in\C^*$).
Then $x^r p(x)^n = \theta\prod_{\alpha}
 (A - \alpha x^s)^{n_\alpha}/x^{s\sum_{\alpha} n_\alpha}$.
Note that each $p_\alpha:=A-\alpha x^s$ is a polynomial, and no
two $p_\alpha$'s have a common root, and $x=0$ is not a root of any
$p_\alpha$.  Thus, for each $\alpha$,
every root of $p_\alpha^{n_\alpha}$ has multiplicity
divisible by $n$, so every
root of $p_\alpha$ has multiplicity divisible by
$n/\gcd(n,n_\alpha)$.

Suppose $\alpha,\beta$ are distinct roots of $g$ such that neither
$n_\alpha$ nor $n_\beta$ is divisible by $n$.  Then
$A-\alpha x^s = a^i$ and $A-\beta x^s=b^j$ where $a,b\in\C[x]$ and
$i,j>1$ are divisors of $n$.  Thus $a^i-b^j=(\beta-\alpha)x^s$,
so $\hat{a}:=a(x^i)/(\beta-\alpha)^{1/i}$ and
$\hat{b}:=b(x^i)/(\beta-\alpha)^{1/j}$ satisfy
$\hat{a}^i-\hat{b}^j=x^{is}$.  Note that $x\nmid\hat{a}\hat{b}$.
Now 
\[\hat{b}^j=\hat{a}^i-(x^s)^i=\prod_{\zeta^i=1}(\hat{a}-\zeta x^s),\]
and the various polynomials $\hat{a}-\zeta x^s$ are coprime
(since $x\nmid\hat{a}$), so for each $\zeta$ we have
$\hat{a}-\zeta x^s = A_{\zeta}^j$ for some $A_{\zeta}\in\C[x]$.
Moreover, we may assume that $\hat{b}=\prod_{\zeta} A_{\zeta}$.
Pick some $\zeta\ne 1$ with $\zeta^i=1$.
Since $x\nmid\hat{a}$, we have $x\nmid A_1 A_{\zeta}$ and
$\gcd(A_1,A_{\zeta})=1$.  But
\[
\prod_{\xi^j=1}(A_1-\xi A_{\zeta}) = A_1^j-A_{\zeta}^j = (\zeta-1) x^s,\]
and any two polynomials $A_1-\xi A_{\zeta}$ are coprime, so every
$A_1-\xi A_{\zeta}$ is an $s\tth$ power.  Since each of these polynomials
divides $x^s$, it follows that one of them is a constant times $x^s$, and
the rest are constants.  But since at least one of $A_1$ and $A_{\zeta}$
is nonconstant, there is at most one $\xi$ for which $A_1-\xi A_{\zeta}$
is constant, whence $j=2$.  Similarly $i=2$, so $\zeta=-1$.
Solving for $A_1$ and $A_{\zeta}$, and then $\hat{a}$ and $\hat{b}$,
we find that $a=\gamma+\delta x^s$ and $b=\pm(\gamma-\delta x^s)$
for some $\gamma,\delta\in\C^*$.
Since $a^2-b^2=(\beta-\alpha)x^s$,
we have $4\gamma\delta=\beta-\alpha$; moreover,
$A=\alpha x^s+a^2=\delta^2 x^{2s} + (\beta+\alpha)x^s/2 + \gamma^2$.
Conversely, given $A$ and $s$, this last equation determines the values
of $\alpha+\beta$, $\gamma^2$, and $\delta^2$, and hence also
$16\gamma^2\delta^2=(\beta-\alpha)^2=(\alpha+\beta)^2-4\alpha\beta$ and finally
$\alpha\beta$.  Thus $A$ and $s$ uniquely determine the set $\{\alpha,\beta\}$.
It follows that $n\mid n_{\chi}$ for every root $\chi$ of $g$ besides
$\alpha$ and $\beta$, whence $g=((x-\alpha)(x-\beta)p^2)^{n/2}$
for some $p\in\C[x]$.  But then $n\mid\deg(g)$, so the order of the
pole of $x^r p^n$ at $x=0$ is divisible by $n$, but this order is $-r$,
contradiction.

This last argument also implies that $g$ is not an $n\tth$ power,
so $g$ has a unique root $\alpha$ for which $n\nmid n_{\alpha}$.
Moreover, for this $\alpha$ we have $\gcd(n,n_{\alpha})=1$.
Thus $g=(x-\alpha)^{n_\alpha} G^n$ for some $G\in\C[x]$,
and $A-\alpha x^s=H^n$ for some $H\in\C[x]$, whence
$h=-\alpha+H^n/x^s$, as desired.
\end{proof}

To determine the decompositions of Laurent polynomials
of the form $x^r p(x^n)$, we use the following result of
Avanzi and Zannier \cite[\S 5]{AZ}:

\begin{prop}[Avanzi--Zannier]\label{AZprop2}
Let $g\in\bC[x]$ be indecomposable, and suppose
$h_1,h_2\in\bC(x)\setminus\bC$ satisfy $g\circ h_1=\gamma g\circ h_2$ where
$\gamma\in\C^*\setminus\{1\}$.
Then $(g,h_1,h_2)=(\theta G\circ\mu,\,\mu^{-1}\circ H_1\circ H,\,
\mu^{-1}\circ H_2\circ H)$ where $\theta\in\C^*$, $\mu\in\bC[x]$ is linear,
$H\in\bC(x)\setminus\bC$, and one of the following occurs:
\begin{enumerate}
\item[(\thethm.1)] $H_2=x$, $\,H_1=\delta x$, and $G\in x^r\C[x^n]$,
where $r\in\bZ_{>0}$, $\,\delta^r=\gamma$, $\,n\in\bZ_{\ge 0}$, and
$\delta^n=1$;
\item[(\thethm.2)] $G=D_n$ with $n$ an odd prime, $\gamma=-1$,
$\,H_1=x+1/x$, and $H_2=H_1\circ\delta x$ where $\delta^n=-1$;
\item[(\thethm.3)] $H_1=(1-\delta x^m)/(\delta x^{m+n}-1)$,
$\,H_2=-1+(x^n-1)/(\delta x^{m+n}-1)$,
and $G=x^m(x+1)^n$,
where $m,n\in\bZ_{>0}$ are coprime and $\delta^n=\gamma$;
\item[(\thethm.4)] $G=D_3(x)+\delta$,
where $\delta\in\C\setminus\{0,2,-2\}$ and either
\begin{enumerate}
\item[(i)]$\gamma=(\delta+2)/(\delta-2)$,
$\,H_1=-1+3(\gamma x^2+1)/(\gamma x^3+1)$, and
$H_2=-2+3(1-x)/(\gamma x^3+1)$; or
\item[(ii)]$\gamma=(\delta-2)/(\delta+2)$,
$\,H_1=-2+3\gamma(1-x)/(x^3+\gamma)$, and
$H_2=-1+3(x^2+\gamma)/(x^3+\gamma)$;
\end{enumerate}
\item[(\thethm.5)] $G=x^4-\frac{4}{3}(\alpha +1)x^3+2\alpha x^2$,
$\,H_1=\frac{(E-\alpha)(E-\frac{1}{\alpha})(x-\frac{6\alpha}{x})+
4(\alpha+1)(E^3+1)}{6(E^4+1)}$, and $H_2=EH_1$,
where $\gamma=-1$, $\,\alpha^4+1=2(\alpha^3+\alpha)$, and 
\[
E=\frac{-1}{2\sqrt{2(2\alpha^2-5\alpha+2)}}\left(x+\frac{6\alpha}x\right)-
\left(\alpha+\frac{1}{\alpha}\right);\]
\item[(\thethm.6)] $G=x^4-\frac{4}{3}(\alpha+\beta)x^3+2\alpha\beta x^2+1$,
where $\omega=e^{2\pi i/3}$, $\,\gamma\in\{\omega,\omega^2\}$,
$\,(\alpha+\omega^2)^3=-2$, and
 $\beta=(1-\alpha)\omega-1$;
if $\gamma=\omega$ then $H_2=\omega^2(H_1-\alpha)E$ and
\[ H_1 = \frac{(E^2+p E+\frac{i}{\sqrt{3}}\alpha^2-w(\alpha-1))U+
\frac{2i}{\sqrt{3}}((\alpha-1)E^3-\omega(\alpha-\omega))}{E^4-1}+\alpha,\]
where $E=(x-\delta/x)/2+p$ and $U=(x+\delta/x)/(2\sqrt{-3(\alpha-1)/2})$
with $p=-\frac{i\omega}{\sqrt{3}}\alpha^2-\omega(\alpha-1)$ and
$\delta=-\omega(\alpha^2-i\sqrt{3}\alpha+3\omega)$;
if $\gamma=\omega^2$ then exchange the above $H_1$ and $H_2$;
\item[(\thethm.7)] $G=x(x+\alpha)^2(x+1)^2$ and $H_2=-Z^2 H_1$,
where $\gamma=-1$ and $Z:=(x-\frac{251+7\xi}{x}+6-2\xi)/32$ with
$\xi^2+\xi+4=0$ and $\alpha^2-\frac{22+5\xi}{9}\alpha+1=0$, and
\[H_1 = \frac{(\alpha+1)(Z^3+1)+(\alpha-1)(Z^2-\xi Z+1)U}{2(Z^5-1)}\]
with $U:=(x+\frac{251+7\xi}{x})/32$.
\end{enumerate}
\end{prop}

\begin{remark}
In the above statement we have implicitly made several corrections
to the results stated in \cite{AZ}.  Specifically, in the definition
of $P_4$ in \cite{AZ}, the equation for $\xi$ should be $\xi^2-2\xi-2=0$.
Our other corrections refer to \cite[Prop.~5.6]{AZ}.  In cases (1) and (3) of
that result, $g_1$ and $h_1$ should be switched; in case (8), $U$
should be replaced by $U/16$; and in case (7), the sign preceding $2/3$
in the expression for $g_1$ should be `$+$', and also an additional
comment must be made for the case $c=\omega^2$.  We also combined
case (1) of \cite[Prop.~5.2, 5.6]{AZ} with case (3), and we combined case (2)
with cases (3) and (4).
\end{remark}

Avanzi and Zannier \cite[Thm.~2]{AZ} generalized Proposition~\ref{AZprop2}
to the case of decomposable $g$, obtaining a recursive description of the
possible polynomials $g$.  In case the genus-zero factor can be
parametrized by Laurent polynomials, we require the following
non-recursive description.

\begin{prop}\label{AZ2me}
Let $g\in\C[x]$ satisfy $\deg(g)>1$, and let $h_1,h_2\in\cL\setminus\bC$ and
$\gamma\in\bC\setminus\{1\}$ satisfy $g\circ h_1=\gamma g\circ h_2$.
Then, after replacing $(g,h_1,h_2)$ by
$(g\circ\mu,\,\mu^{-1}\circ h_1\circ\theta x,\,\mu^{-1}\circ h_2\circ\theta x)$
for some $\theta\in\C^*$ and some linear $\mu\in\C[x]$,
one of the following holds (where $n\in\bZ_{\ge 0}$ and $r,m\in\bZ_{>0}$):
\begin{enumerate}
\item[(\thethm.1)] $h_1=\alpha h_2$ and $g\in x^r\C[x^n]$, where $\alpha^n=1$
and $\alpha^r=\gamma$;
\item[(\thethm.2)] $h_1=x^m+1/x^m$,
$\,h_2=h_1\circ \alpha x$, and
$g=G\circ D_n$, where $\gamma=-1$, $\,G\in x\C[x^2]$, and
$\alpha^{nm}=-1$;
\item[(\thethm.3)] $h_1=x^m+1/x^m$,
$\,h_2=(x^m-1/x^m)/\sqrt{\alpha}$,
and $g=G\circ (\frac{(1-\alpha)x^2}{2}-2)$,
where $G\in x^r\C[x^n]$, $\,\alpha^r=\gamma$, and
$\alpha^n=1$ but $\alpha\ne -1$.
\end{enumerate}
\end{prop}

\begin{proof}
Write $g=g_1\circ \dots\circ g_s$ where the $g_i$ are indecomposable
polynomials.  Let $j$ be the largest integer $\le s$ for which
$H_1:=g_{j+1}\circ\dots\circ g_s\circ h_1$ and
$H_2:=g_{j+1}\circ\dots\circ g_s\circ h_2$ satisfy
$g_j\circ H_1=\nu\circ g_j\circ H_2$ for some linear
$\nu\in\C[x]$, and put
$G=g_1\circ\dots\circ g_{j-1}$.  Writing $\nu(x)=\alpha x+\beta$
and comparing leading coefficients in the identity
$G\circ\nu=\gamma G$,
we see that $\alpha^{\deg(G)}=\gamma\ne 1$, so $\alpha\ne 1$.
Now put $\lambda:=x+\beta/(\alpha-1)$, so
$\lambda\circ\nu=\alpha\lambda$; replacing $G$ and $g_j$ by
$G\circ\lambda^{-1}$ and $\lambda\circ g_j$, we have
$g_j\circ H_1=\alpha g_j\circ H_2$, so
$G(\alpha x)=\gamma G(x)$.
Hence $G\in x^r\C[x^n]$ for some $r>0$ and $n\ge 0$ such that
$\alpha^n=1$ and $\alpha^r=\gamma$.
If $h_1=\hat\nu\circ h_2$ with $\hat\nu\in\C[x]$ linear, then this argument
shows that (\ref{AZ2me}.1) holds.  Henceforth assume
there is no such $\hat\nu$, so there is no linear $\hat\nu\in\C[x]$ such that
$H_1=\hat\nu\circ H_2$.

By Proposition~\ref{AZprop2}, there exist $\hat\theta\in\C^*$,
$\,H\in\bC(x)\setminus\bC$, and a linear $\mu\bC[x]$ such that
\begin{align*}
g_j&=\hat\theta \hat{g_j}\circ\hat\mu \\
H_1&= \hat\mu^{-1}\circ \hat{H_1}\circ H \\
H_2&= \hat\mu^{-1}\circ \hat{H_2}\circ H,
\end{align*}
where $\hat{g_j}$, $\hat{H_1}$, and $\hat{H_2}$ satisfy the conditions
required of $G$, $H_1$, and $H_2$ in one of
(\ref{AZprop2}.1)--(\ref{AZprop2}.7).
By replacing $G$ by $G\circ \hat\theta x$, we may replace $g_j$ by $\hat{g_j}$
while also replacing $H_1$ and $H_2$ by $\hat{H_1}\circ H$ and
$\hat{H_2}\circ H$.

Since $H_1,H_2\in\cL$ have at most two poles,
also $\hat{H_1}$ and $\hat{H_2}$ have at most two poles.  This rules out
(\ref{AZprop2}.4)--(\ref{AZprop2}.7).  In (\ref{AZprop2}.3) it implies
$m=n=1$, so $g_j=x^2+x$, $\,\hat{H_1}=(1-\alpha x)/(\alpha x^2-1)$
and $\hat{H_2}=(x-\alpha x^2)/(\alpha x^2-1)$.  Putting
\[
\mu_1=\frac{4x+2}{\sqrt{1-\alpha}} \quad\text{and}\quad
\mu_2=\frac{1}{\sqrt{\alpha}}\frac{x(1+\sqrt{\alpha})+\sqrt{1-\alpha}}
  {x(1+\sqrt{\alpha})-\sqrt{1-\alpha}},\]
we have
\begin{align*}
8g_j\circ\mu_1^{-1} &= \frac{1-\alpha}{2}x^2-2\\
\mu_1\circ\hat{H_1}\circ\mu_2 &= x+\frac{1}{x}\\
\mu_1\circ\hat{H_2}\circ\mu_2 &=
 \frac{1}{\sqrt{\alpha}}\left(x-\frac{1}{x}\right).
\end{align*}
Now replace $G$ by $G\circ 8x$ and $g_j$ by $g_j\circ\mu_1^{-1}$,
while also replacing $\hat{H_1}$ and $\hat{H_2}$ by
$\mu_1\circ\hat{H_1}\circ\mu_2$ and $\mu_2\circ\hat{H_2}\circ\mu_2$
(and replacing $H$ by $\mu_2^{-1}\circ H$).
Thus we have $g_j=\frac{1-\alpha}{2}x^2-2$,
$\,\hat{H_1}=x+x^{-1}$, and $\hat{H_2}=(x-x^{-1})/\sqrt{\alpha}$.
Since $H_1=\hat{H_1}\circ H$ has no poles besides $0$ and $\infty$,
and $\hat{H_1}$ has poles at $0$ and $\infty$, the full $H$-preimage
of $\{0,\infty\}$ is $\{0,\infty\}$, so $H=(\theta x)^m$ for some
nonzero $m\in\bZ$ and $\theta\in\C^*$.  If $m<0$ then replace
$H$ by $(\theta x)^{-m}$ and $\hat{H_2}$ by $-\hat{H_2}$, thereby
preserving the compositions $\hat{H_1}\circ H$ and
$\hat{H_2}\circ H$.  Thus we may assume $m>0$ by making the
appropriate choice of $\sqrt{\alpha}$.  Now
$H_1=\hat{H_1}\circ H=(x^m+x^{-m})\circ\theta x$ and
$H_2=(x^m-x^{-m})/\sqrt{\alpha}\circ\theta x$.
Write $R=g_{j+1}\circ\dots\circ g_s$, so
\begin{align*}
R\circ h_1 &= H_1 = \left(x^m+\frac{1}{x^m}\right)\circ \theta x \\
R\circ h_2 &= H_2 = \frac{1}{\sqrt{\alpha}}\left(x^m-\frac{1}{x^m}\right)
                          \circ \theta x.
\end{align*}
By Lemma~\ref{Dickson}, we have
$R=D_{m/d}\circ\mu$ where $d\mid m$ and
$\mu\in\C[x]$ is linear; moreover, $h_1=\mu^{-1}\circ (x^d+1/x^d)\circ 
\theta x$.  Since $R\circ h_2 = \frac{x^m+x^{-m}}{i\sqrt{\alpha}}\circ
\theta i^{1/m} x$, Lemma~\ref{Dickson} implies that
$R=D_{m/d}(x)/(i\sqrt{\alpha})\circ \tilde\mu$ and
$h_2=\tilde\mu^{-1}\circ (x^d+x^{-d})\circ\theta i^{1/m} x$
for some linear $\tilde\mu\in\C[x]$.
Equating coefficients in the identity
$D_{m/d}\circ\mu = R = D_{m/d}/(i\sqrt{\alpha})\circ\tilde\mu$,
we see that either $\alpha=-1$ or $m=d$.  If $m=d$ then
$g=G\circ g_j\circ\mu$ and $h_1=\mu^{-1}\circ (x^m+x^{-m})\circ\theta x$
and $h_2=\mu^{-1}\circ (x^m-x^{-m})/\sqrt{\alpha}\circ\theta x$,
as in (\ref{AZ2me}.3).  Now assume $m\ne d$, so $\alpha=-1$,
whence $g_j=D_2$.  Replacing $g$, $h_1$ and $h_2$ by
$g\circ\mu^{-1}$, $\,\mu\circ h_1\circ x/\theta$, and
$\mu\circ h_2\circ x/\theta$, we have $g=G\circ D_{2m/d}$ and
$h_1=x^d+x^{-d}$ and $h_2=\pm h_1\circ i^{1/m} x$.  Thus
$h_2=h_1\circ\hat\alpha x$ where $\hat\alpha^{2m}=-1$, and we
have obtained (\ref{AZ2me}.2) with $n=2$.

Now assume $g_j$, $\hat{H_1}$ and $\hat{H_2}$ satisfy
(\ref{AZprop2}.1).  Then $\hat{H_1}=\delta x$ and $\hat{H_2}=x$
for some $\delta\in\C^*$, so $H_1=\hat{H_1}\circ H=\delta H_2$,
contradicting our hypothesis to the contrary.

Finally, assume $g_j$, $\hat{H_1}$ and $\hat{H_2}$ satisfy
(\ref{AZprop2}.3).  Thus $\alpha=-1$ and $g_j=D_p$ with $p$ an
odd prime, and moreover $\hat{H_1}=x+1/x$ and
$\hat{H_2}=\hat{H_1}\circ\delta x$ where $\delta^p=-1$.
Since $H_1=\hat{H_1}\circ H$ is a Laurent polynomial, we must have
$H=(\theta x)^m$
for some nonzero $m\in\bZ$ and $\theta\in\C^*$.  
If $m<0$ then we can replace $m$ by $-m$ if we replace
$\delta$ and $\theta$ by $1/\delta$ and $1/\theta$;
since these changes do not affect $H_1$ or $H_2$, we may assume $m>0$.
Write
$R=g_{j+1}\circ\dots\circ g_s$, so
$R\circ h_1=(x+1/x)\circ (\theta x)^m$
and $R\circ h_2=(x+1/x)\circ\delta (\theta x)^m$.
By Lemma~\ref{Dickson}, we have
$R=D_{m/d}\circ\mu$ where $d\mid m$ and $\mu\in\C[x]$ is linear;
moreover, $h_1=\mu^{-1}\circ (x^d+1/x^d)\circ \theta x$.
Likewise $R=D_{m/d}\circ\tilde\mu$ for some linear $\tilde\mu\in\C[x]$,
and moreover
$h_2=\tilde\mu^{-1}\circ (x^d+1/x^d)\circ x\theta\delta^{1/m}$.
The identity
$D_{m/d}\circ\mu=R=D_{m/d}\circ\tilde\mu$ implies that
$\tilde\mu=\epsilon\mu$ with $\epsilon\in\{1,-1\}$ and
$\epsilon^{m/d}=1$.  After replacing $g$, $h_1$ and $h_2$ by
$g\circ\mu^{-1}$, $\,\mu\circ h_1\circ x/\theta$, and
$\,\mu\circ h_2\circ x/\theta$, we have
$g=G\circ D_{pm/d}$ and $h_1=x^d+1/x^d$ and
$h_2=\epsilon h_1\circ x\delta^{1/m}$, so
$h_2=h_1\circ \hat\alpha x$ where $\hat\alpha^{mp}=-1$.
Thus we have (\ref{AZ2me}.2).
\end{proof}

We can now describe the decompositions of Laurent polynomials of
the form $x^r p(x^n)$:

\begin{prop}\label{Ritttwisthard}
Let $n,r\in\bZ$ satisfy $n>1$ and $n\nmid r$, and pick $p\in\C[x]$
with $x\nmid p$.
Suppose $g,h\in\C(x)$ satisfy $g\circ h = x^r p(x^n)$.  Then there
is a degree-one $\mu\in\C(x)$ such
that, after replacing $g$ and $h$ by $g\circ\mu$ and $\mu^{-1}\circ h$,
one of the following occurs (with $s,t,m\in\bZ$ and $m>0$):
\begin{enumerate}
\item[(\thethm.1)] $g\in x^s\C[x^m]$ and $h\in x^t\C[x^n]$ where
   $n\mid mt$;
\item[(\thethm.2)] $g=G\circ D_t$ and $h=(x^m+1/x^m)\circ\theta x$ where
$G\in x\C[x^2]$ and $mt\equiv r\equiv n/2\pmod{n}$, with $n$ even, $t>0$,
and $\theta\in\C^*$.
\end{enumerate}
Moreover, if $g\in\C[x]$ and $h\in\cL$ then we may choose
$mu\in\C[x]$.
\end{prop}

\begin{proof}
By Lemma~\ref{basic}, we may assume $g,h\in\cL$ and either $g\in\C[x]$
or $h=x^t$ with $t\in\bZ_{>0}$.  In the latter case the condition
$x^r p(x^n)\in\C[x^t]$ implies $t\mid r$ and $p=P(x^{t/\gcd(n,t)})$
with $P\in\C[x]$.  Thus $g=x^{r/t} P(x^{n/\gcd(n,t)})$, as in
(\ref{Ritttwisthard}.1).  Henceforth assume $g\in\C[x]$.  If
$\deg(g)=1$ then we may assume $g=x$, so again (\ref{Ritttwisthard}.1)
holds.  Now assume $\deg(g)>1$.

Let $\zeta$ be a primitive $n\tth$ root of unity.
Then $g\circ h(\zeta x) = \zeta^r g\circ h(x)$, and $\gamma:=\zeta^r\ne 1$.
Write $h_2:=h(x)$ and $h_1:=h(\zeta x)$, so $g\circ h_1= \gamma g\circ h_2$.
By Proposition~\ref{AZ2me}, there exist $\theta\in\C^*$ and a
linear $\mu\in\C[x]$ such that, after replacing $g,h_1,h_2$ by
$g\circ\mu$, $\,\mu^{-1}\circ h_1\circ\theta x$, and
$\mu^{-1}\circ h_2\circ\theta x$,
one of (\ref{AZ2me}.1)--(\ref{AZ2me}.3) holds.
We will use the equation $h_1=h_2\circ \zeta x$ to analyze these
possibilities.

If (\ref{AZ2me}.1) holds then $\alpha h_2=h_1=h_2\circ \zeta x$, so
$h_2\in x^t\C[x^n]$ with $\zeta^t=\alpha$;
here also $g\in x^s\C[x^m]$ where $\alpha^m=1$ and $\alpha^s=\gamma$.
Thus $\zeta^{tm}=1$, so we have (\ref{Ritttwisthard}.1).

If (\ref{AZ2me}.2) holds then
\[x^m+\frac{1}{x^m} = h_1 = h_2\circ \zeta x =
\left(x^m+\frac{1}{x^m}\right)\circ \alpha\zeta x,\]
so $(\alpha\zeta)^m=1$.
Here $g=G\circ D_t$ where $\gamma=-1$ and $G\in x\C[x^2]$,
and $\alpha^{mt}=-1$.  Thus $\zeta^{mt}=-1$, and we have
(\ref{Ritttwisthard}.2).

If (\ref{AZ2me}.3) holds then, for some $\alpha\ne -1$, we have
\[x^m+\frac{1}{x^m} = h_1 = h_2\circ \zeta x = 
\frac{1}{\sqrt{\alpha}}\left(x^m-\frac{1}{x^m}\right)\circ \zeta x,
\] so $\zeta^m=\sqrt{\alpha}=-1/\zeta^m$.  But then
$\alpha=\zeta^{2m}=-1$, contradiction.
\end{proof}

Next we consider decompositions of $(x^2-4)p(x)^2$ with $p\in\C[x]$;
since these are polynomials, Ritt's results provide information
about their decompositions, but we go further by precisely describing
the shape of every decomposition:

\begin{prop}\label{decompweird1}
Let $g,h,p\in\C[x]\setminus\{0\}$ satisfy $g\circ h = (x^2-4)p(x)^2$.
Then, after replacing $g$ and $h$ by $g\circ\mu$ and $\mu^{-1}\circ h$
for some linear $\mu\in\C[x]$, there exist $B,D\in\C[x]$ and $n\in\bZ_{>0}$ such
that one of the following holds:
\begin{enumerate}
\item[(\thethm.1)] $g=xB^2$ and $h=(x^2-4)D^2$;
\item[(\thethm.2)] $g = (x^2-4)B^2$ and $h=D_n$.
\end{enumerate}
\end{prop}

\begin{remark} To verify that the polynomials $g$ and $h$ in
(\ref{decompweird1}.2) satisfy $g\circ h = (x^2-4)p(x)^2$ for suitable $p$,
note that $D_n^2-4=(x^2-4)E_{n-1}^2$, where the polynomial $E_{n-1}$ is a
`Dickson polynomial of the second kind', and is defined by the functional
equation
$E_{n-1}(x+x^{-1}) = (x^n-x^{-n})/(x-x^{-1})$.
\end{remark}

\begin{proof}[Proof of Proposition~\ref{decompweird1}]
Write
$g=AB^2$ and $h=CD^2$ with $A,B,C,D\in\C[x]$ and $A,C$ squarefree and monic.
Then $(x^2-4)p(x)^2 = A(h)\cdot B(h)^2$, so $A(h)$ is a square
times $x^2-4$.  Write $A(x)=\prod_\alpha (x-\alpha)$, where the product
ranges over the roots of $A$, and
write $h-\alpha = E_{\alpha}^2 F_{\alpha}$ with $E_\alpha,F_\alpha\in
\C[x]$ and $F_\alpha$ squarefree and monic.  
For distinct roots $\alpha,\alpha'$ of $A$, plainly
$h-\alpha$ and $h-\alpha'$ are coprime, so
$\gcd(E_\alpha,E_{\alpha'})=1=\gcd(F_\alpha,F_{\alpha'})$.
Since $A(h)=\prod_\alpha E_{\alpha}^2F_{\alpha}$ is a square times $x^2-4$,
and the various polynomials $F_{\alpha}$ are monic, squarefree and coprime,
we have $x^2-4=\prod_{\alpha} F_{\alpha}$.
Moreover, differentiating the equation $h-\alpha=E_{\alpha}^2 F_{\alpha}$
implies $E_{\alpha}\mid h'$, and since the various polynomials $E_\alpha$
are coprime, we have $\prod_\alpha E_\alpha \mid h'$.
Writing $n=\deg(h)$ and $r=\deg(A)$,
it follows that $n-1\ge\sum_\alpha\deg(E_\alpha)$, and since
\[
nr=\deg(h\circ A)=\sum_\alpha\deg(h-\alpha)=
2+2\sum_\alpha\deg(E_\alpha),\]
we conclude that $r\le 2$.
If $r=1$ then, after replacing $g$ and $h$ by $g\circ\mu$ and
$\mu^{-1}\circ h$ for a suitable linear $\mu\in\C[x]$, we may assume
$A=x$; but then $C=F_0=x^2-4$, so we have (\ref{decompweird1}.1).
Now assume $r=2$, so, after inserting a
linear and its inverse between $g$ and $h$ as above, we may assume $A=x^2-4$.
There are four possibilities:
\begin{enumerate}
\item[(i)] $F_2=x^2-4$ and $F_{-2}=1$;
\item[(ii)] $F_2=x-2$ and $F_{-2}=x+2$;
\item[(iii)] $F_2=x+2$ and $F_{-2}=x-2$; or
\item[(iv)] $F_2=1$ and $F_{-2}=x^2-4$.
\end{enumerate}
By replacing $g$ and $h$ by $g\circ (-x)$ and $(-x)\circ h$,
we may assume that (i) or (ii) holds.  In either case, the cover
$h:\Line\to\Line$ is totally ramified over $\infty$, and every point
lying over $2$ or $-2$ has even ramification index except $2$ and $-2$.
This data determines $h$ up to composition on both sides with linears,
as was first shown by Ritt \cite{Ritt}, and as has been reproved in every
proof of Ritt's results.  Thus,
$h=\nu_1\circ D_n \circ \nu_2$ for some linear $\nu_1,\nu_2\in\C[x]$.
In case (ii) we have $h-2=(x-2)E_2^2$, so $n$ is odd;
if $n=1$ then $E_2$ is a constant, and since $(x+2)\mid(h+2)$ we must
have $E_2=\pm 1$, so $h=x$ and (\ref{decompweird1}.2) holds.
If (ii) holds with $n>1$ then $2$ and $-2$ are the unique finite
branch points of $h:\Line\to\Line$, and their unique unramified preimages
are $2$ and $-2$, respectively.  Since $D_n$ has the same property,
each $\nu_i$ preserves $\{2,-2\}$, hence equals $\pm x$,
and we must have $\nu_2=\nu_1$.  Since
$-D_n(-x)=D_n(x)$, this gives (\ref{decompweird1}.2).  In case (i), $n$ is
even; if $n=2$ then $-2$ is the unique finite branch point of both $h$ and
$D_n$, so $\nu_1$ fixes $-2$ and thus $\nu_1=-2+\beta\cdot (x+2)$.
Since $h(\pm 2)=2$ and $D_2=x^2-2$, we find that $\nu_2=\alpha x$
where $\beta=1/\alpha^2$, which implies $h=D_2$ as desired.  Now suppose
$n>2$.  Then both $h$ and $D_n$ have $2$ and $-2$ as their unique finite
branch points, and all of their preimages are ramified except for $\pm 2$,
both of which lie over $2$.  Thus $\nu_1$ fixes $2$ and $-2$, so $\nu_1=x$.
Also $\nu_2$ preserves $\{2,-2\}$, so $\nu_2=\pm x$, whence $h=D_n$.
\end{proof}

Finally, we determine the decompositions of the other Laurent polynomials
in (\ref{Lbidec}.2), namely 
$(x-1/x)\cdot p(x+1/x)$ with $p\in\C[x]$.  As we noted in (\ref{q}),
composition with a degree-one rational function transforms these into
the form $x q(x^2)$, but the resulting $q\in\C(x)$ is not a Laurent
polynomial.

\begin{prop}\label{decompweird2}
Let $g,p\in\C[x]$ and $h\in\cL$ satisfy $p\ne 0$ and
$g\circ h = (x-1/x)\cdot p(x+1/x)$.  Then there exist $\mu,q\in\C[x]$
with $\mu$ linear such that one of the following holds:
\begin{enumerate}
\item[(\thethm.1)] $\mu^{-1}\circ h=(x-1/x)\cdot q(x+1/x)$ and
$g\circ\mu\in x\C[x^2]$ is an odd polynomial;
\item[(\thethm.2)] $\mu^{-1}\circ h=\frac{x^m}{\sqrt{\gamma}}+
\frac{\sqrt{\gamma}}{x^m}$ and $g\circ\mu=G\circ D_n$ with
$G\in x\C[x^2]$ and $\gamma^n=-1$.
\end{enumerate}
\end{prop}

\begin{remark}
We note that the examples in (\ref{decompweird2}.2) do satisfy the hypotheses:
for, $f:=g\circ h = G\circ D_n\circ (x+1/x)\circ x^m/\sqrt{\gamma}$.
Writing $I=\sqrt{\gamma}^n$, we have $I^2=-1$, so
$f = G\circ (x+1/x)\circ I x^{nm} =
G\circ I(x-1/x)\circ x^{nm}$.  There is a polynomial $E_{nm-1}$ (the Dickson
polynomial of the second kind) satisfying $(x-1/x)\circ x^{nm} =
(x-1/x)E_{nm-1}(x+1/x)$.  Since $G$ is odd, it follows that
$f(x)=(x-1/x)\cdot p(x+1/x)$ for some $p\in\C[x]$.
\end{remark}

\begin{proof}[Proof of Proposition~\ref{decompweird2}]
Write $f=(x-1/x)\cdot p(x+1/x)$.  Since $f(1/x)=-f(x)$ (and $f\in\cL$),
we can write
$f(x)=F(x)-F(1/x)$ with $F\in x\C[x]$.
Write the leading terms of $F$ and $g$ as $\beta x^s$ and $\theta x^r$.
Viewing $f$ as a finite Laurent series, its
highest and lowest-degree terms have degrees $s$ and
$-s$, so we can write
$h=\delta(x^e+\delta_1 x^{e-1} + \dots + \delta_{e-1} x) + \xi +
   \zeta(x^{-e} + \zeta_1 x^{1-e} + \dots + \zeta_{e-1} x^{-1})$
with $\delta,\zeta\in\C^*$ and $\delta_i,\zeta_i,\xi\in\C$, where $e=s/r$.
Then $\delta^r=\beta/\theta=-\zeta^r$, and moreover the $\delta_i$ are
uniquely determined by $F$, since the coefficients of $x^{s-1},\dots,x^{s-e+1}$
in the congruence
$(x^e+\delta_1 x^{e-1}+\dots+\delta_{e-1} x)^r\equiv F/\beta \pmod{x^{s-e}}$
successively determine $\delta_1,\dots,\delta_{e-1}$.
Since the $\zeta_i$ are determined by the same congruence, we have
$\zeta_i=\delta_i$, whence $h=H(x)+\gamma H(1/x)+\xi$ with $H\in x\C[x]$ and
$\gamma^r=-1$.
Since $f(1/x)=-f(x)$, we have $g\circ h(x)=-g\circ h(1/x)$.
By Proposition~\ref{AZ2me}, there exist $\hat\theta\in\C^*$ and a linear
$\mu\in\C[x]$ such that one of (\ref{AZ2me}.1)--(\ref{AZ2me}.3) holds for
$\hat g:=g\circ\mu$, 
$h_1:=\mu^{-1}\circ h\circ\hat\theta x$, and
$h_2:=\mu^{-1}\circ h\circ (\hat\theta x)^{-1}$.
Write $\hat H(x) = \mu^{-1}\circ H(x) - \mu^{-1}(0)$, so $\hat H\in x\C[x]$
and $h_1=\hat H(\hat\theta x) + \gamma \hat H(1/(\hat\theta x)) +
\mu^{-1}(\xi)$ and
$h_2=\hat H(1/(\hat\theta x)) + \gamma \hat H(\hat\theta x)+\mu^{-1}(\xi)$.

In case (\ref{AZ2me}.1) we have $h_1=\alpha h_2$, where $\alpha\ne 1$.
Comparing the terms of highest and lowest degrees in this identity
gives $\frac{1}{\alpha}\cdot\mu^{-1} = \mu^{-1}\circ\gamma x = \alpha\cdot
\mu^{-1}$, so $\alpha=\gamma=-1$.  Now (\ref{AZ2me}.1) implies $\hat g\in
x\C[x^2]$.  Since $h_1$ and $h_2=-h_1$ both have
constant term $\mu^{-1}(\xi)$, this term must be zero, so 
$h_1(x)=\hat{H}(x)-\hat{H}(1/x)$.  Letting $\sigma$ be the automorphism
of $\C(x)$ mapping $x\mapsto 1/x$, we see that $R:=h_1(x)/(x-1/x)$ is 
fixed by $\sigma$, and thus lies in the fixed field $\C(x)^\sigma=
\C(x+1/x)$.  Thus $R=q(x+1/x)$ for some $q\in\C(x)$.
The only poles of $1/(x-1/x)$ are $1$ and $-1$, both of which have
order $1$; since $h_1(1)=h_1(-1)=0$, neither $1$ nor $-1$ is a pole of $R$,
so $R$ has no poles besides $0$ and $\infty$.  Since $R=q(x+1/x)$, and
the images of $0$ and $\infty$ under $x+1/x$ are both $\infty$, it follows
that $q$ has no poles besides $\infty$, so $q\in\C[x]$.  This proves
that (\ref{decompweird2}.1) holds.

In cases (\ref{AZ2me}.2) and (\ref{AZ2me}.3) we have
$h_1 = x^m + 1/x^m$, so $x^m=\hat H(\hat\theta x)$ and
$1/x^m=\gamma\hat H(\frac{1}{\hat\theta x})$, whence
$\hat H(x)=(x/\hat\theta)^m$ and $\hat\theta^{2m}=\gamma$.
Thus $h_2 = x^m/\gamma + \gamma/x^m$, which is incompatible with
(\ref{AZ2me}.3), so (\ref{AZ2me}.2) holds.  Moreover, in
(\ref{AZ2me}.2) we must have $\alpha^m=1/\gamma$,
and $\hat g=G\circ D_n$ where $G$ is
an odd polynomial and $\alpha^{mn}=-1$.
This yields (\ref{decompweird2}.2).
\end{proof}


\section{Laurent polynomials with two Type 1 decompositions}
\label{sec-1}

In this section we describe all instances of Laurent polynomials with
two Type 1 decompositions.
Our proofs make crucial use of a result of Bilu and Tichy \cite[Thm.~9.3]{BT},
whose proof relies on Ritt's results among other things.  The statement
of this result involves the general degree-$n$ Dickson polynomial
$D_n(x,\alpha)$ (with $\alpha\in\C$), which is defined by the functional
equation $D_n(z+\alpha/z,\alpha)=z^n+(\alpha/z)^n$ (in this notation,
our previously defined $D_n(x)$ is $D_n(x,1)$).

\begin{prop}[Bilu--Tichy]
\label{BTthm}
Let $g_1,g_2\in \C[x]\setminus\C$, and let $E(x,y)\in\C[x,y]$ be a factor of
$g_1(x)-g_2(y)$.
Suppose that $E(x,y)=0$ is an irreducible curve of genus $0$
which has at most two closed points lying over $x=\infty$.  Then
$g_1=G\circ G_1\circ\mu_1$ and $g_2=G\circ G_2\circ\mu_2$, where
$G,\mu_1,\mu_2\in \C[x]$ with $\mu_1,\mu_2$ linear, and
where either $(G_1,G_2)$ or $(G_2,G_1)$ is in the following list
(in which $p\in\C[x]$ is nonzero, $m,n$ are coprime positive integers,
 and $\alpha,\beta\in\C^*$):
\begin{enumerate}
\item[(\thethm.1)] $(x^n, \,\alpha x^r p(x)^n)$ where $0\le r<n$ and
$\gcd(r,n)=1$;
\item[(\thethm.2)] $(x^2, \,(\alpha x^2+\beta)p(x)^2)$;
\item[(\thethm.3)] $(D_m(x,\alpha^n), \,D_n(x,\alpha^m))$;
\item[(\thethm.4)] $(\alpha^{-m}D_{2m}(x,\alpha),
 \,-\beta^{-n}D_{2n}(x,\beta))$;
\item[(\thethm.5)] $((\alpha x^2-1)^3, \,3x^4-4x^3)$;
\item[(\thethm.6)] $(D_{dm}(x,\alpha^n),
\,-D_{dn}(x\cos(\pi/d),\alpha^m))$ where $d\ge 3$.
\end{enumerate}
Moreover, there exists $(G_1,G_2)$ as above such that $E(x,y)$ is a
factor of $G_1\circ\mu_1(x)-G_2\circ\mu_2(y)$, and such that in all
but the last case $E(x,y)$ is a constant times
$G_1\circ\mu_1(x)-G_2\circ\mu_2(y)$.
\end{prop}

\begin{remark}
In the above result we have corrected an error from \cite{BT},
namely that $a^{n/d}$ and $a^{m/d}$ should be switched in the definition
of `specific pairs' in \cite{BT} in order to make
\cite[Thm.~9.3]{BT} be true.
\end{remark}

Actually Bilu and Tichy proved a version of this result for polynomials
over an arbitrary field of characteristic zero; since we have restricted
to the complex numbers, we can simplify the statement somewhat:

\begin{cor}
\label{BTcor}
Proposition~\ref{BTthm} remains true if we replace
\emph{(\ref{BTthm}.1)}--\emph{(\ref{BTthm}.6)} by the following
(where $m,n\in\bZ_{>0}$ are coprime, and $p\in\C[x]$ is nonzero):
\begin{enumerate}
\item[(\thethm.1)] $(x^n, \,x^r p(x)^n)$ where $0\le r<n$ and $\gcd(r,n)=1$;
\item[(\thethm.2)] $(x^2, \,(x^2-4)p(x)^2)$;
\item[(\thethm.3)] $(D_m(x), \,D_n(x))$;
\item[(\thethm.4)] $((x^2/3-1)^3, \,3x^4-4x^3)$;
\item[(\thethm.5)] $(D_{dm}(x), \,-D_{dn}(x))$ where $d>1$.
\end{enumerate}
\end{cor}

Before proving Corollary~\ref{BTcor}, we recall some basic properties
of Dickson polynomials.  These follow readily from
the definition; for details, and further results, see \cite{ACZ,LMT}.
\stepcounter{lemma}
\begin{gather}
\tag{\thethm.1}\label{D12}D_1(x,\alpha)=x;\quad D_2(x,\alpha)=x^2-2\alpha; \\
\tag{\thethm.2}\label{Dmn}D_{mn}(x,\alpha)=D_m(D_n(x,\alpha),\alpha^n); \\
\tag{\thethm.3}\label{cD}\beta^n D_n(x,\alpha)=D_n(\beta x,\beta^2\alpha).
\end{gather}

\begin{proof}[Proof of Corollary~\ref{BTcor}]
If (\ref{BTthm}.1) holds then (\ref{BTcor}.1) holds, since
$\alpha x^r p(x)^n = x^r (\sqrt[n]{\alpha} p(x))^n$.
Likewise, if (\ref{BTthm}.2) holds then so does (\ref{BTcor}.2) (perhaps
after changing $p$ and $\mu_i$),
since $(\alpha x^2+\beta)p(x)^2 = (x^2-4)\hat{p}(x)^2\circ \gamma x$
where $\gamma^2=-4\alpha/\beta$ and $\hat{p}(x)=(\sqrt{-\beta}/2)
p(x/\gamma)$.
We pass from (\ref{BTthm}.5) to (\ref{BTcor}.4) in a similar manner,
since $(\alpha x^2-1)^3 = (x^2/3-1)^3 \circ \sqrt{3\alpha}x$.
If (\ref{BTthm}.4) holds, we use
(\ref{cD}) with $\gamma^2=1/\alpha$ and $\delta^2=1/\beta$,
getting $\alpha^{-m}D_{2m}(x,\alpha)=D_{2m}(x\gamma)$ and
$-\beta^{-n}D_{2n}(x,\beta)=-D_{2n}(x\delta)$, which yields (\ref{BTcor}.5)
(with $d=2$).

If (\ref{BTthm}.3) holds, let $\gamma$ be a square root of $\alpha$, so
(\ref{cD}) implies $D_m(x,\alpha^n)=\gamma^{nm}D_m(x/\gamma^n)$, whence
$G\circ D_m(x,\alpha^n) = G(\gamma^{nm}x)\circ D_m(x)\circ x/\gamma^n$.
Since we could do the same thing after exchanging $n$ and $m$,
and since this change would not affect $G(\gamma^{nm}x)$, it follows that
(\ref{BTcor}.3) holds here.

If (\ref{BTthm}.6) holds, we again let $\gamma$ be a square root of $\alpha$,
so (\ref{cD}) implies that 
$-D_{dn}(x\cos(\pi/d),\alpha^m)=-\gamma^{dmn}D_{dn}(x\cos(\pi/d)/\gamma^m)$
and $D_{dm}(x,\alpha^n)=\gamma^{dmn}D_{dm}(x/\gamma^n)$.
Thus, after replacing $G(x)$ by $G(\gamma^{dmn}x)$, and composing
$\mu_1$ and $\mu_2$ with $x\cos(\pi/d)/\gamma^m$ and $x/\gamma^n$,
we obtain (\ref{BTcor}.5).
\end{proof}

To describe the Laurent polynomials with two Type 1 decompositions, we
need two more auxiliary results.
The first is a neat observation of Fried's about
factorizations of polynomials of the form
$g_1(x)-g_2(y)$ \cite[Prop.~2]{Friedfac};
we state the refined version given in \cite[Thm.~8.1]{BT}:

\begin{prop}\label{Fried}
For any $G_1,G_2\in\C[x]\setminus\C$, there exist
$a_1,a_2,b_2,b_2\in\C[x]$ such that
\begin{enumerate}
\item[(\thethm.1)] $G_1=a_1\circ b_1$\quad\text{and}\quad $G_2=a_2\circ b_2$;
\item[(\thethm.2)] the splitting field of $a_1(x)-z$ over $\C(z)$ equals
                   the splitting field of $a_2(x)-z$ over $\C(z)$;
\item[(\thethm.3)] the irreducible factors of $G_1(x)-G_2(y)$ are precisely
                   the polynomials $A(b_1(x),b_2(y))$, where $A$ is
                   an irreducible factor of $a_1(x)-a_2(y)$.
\end{enumerate}
\end{prop}

We also require the factorization of $D_n(x)+D_n(y)$; as noted by
Bilu \cite[Prop.~3.1]{Bilu}, (\ref{D12}) and (\ref{Dmn}) imply
$D_{2n}=D_n^2-2$, so for $F_n:=D_n(x)-D_n(y)$ we have
$D_n(x)+D_n(y)=F_{2n}/F_n$, and hence
it suffices to factor $F_n$.  This last factorization is
well-known; for a simple derivation see \cite[Thm.~7]{BhZ}.

\begin{prop}\label{Dfac}
Put
\[
\Phi_n(x,y) = \prod_{\substack{{1\le k<n} \\ {k\equiv 1\bmod 2}}}
(x^2 - xy\cdot 2\cos(\pi k/n) + y^2 - 4\sin^2(\pi k/n)).\]
Then
\[
D_n(x) + D_n(y) =
\begin{cases}
\Phi_n(x,y) & \text{if $n$ is even} \\
(x+y)\Phi_n(x,y) & \text{if $n$ is odd.}
\end{cases}
\]
\end{prop}

We now classify Laurent polynomials with two Type 1 decompositions.

\begin{thm} \label{Type1}
Let $g_1,g_2\in\C[x]\setminus\C$ and $h_1,h_2\in\cL\setminus\C$ satisfy
$g_1\circ h_1=g_2\circ h_2$.  Then, perhaps after switching
$(g_1,g_2)$ and $(h_1,h_2)$, we have
\begin{align*}
g_1&=G\circ G_1\circ\mu_1\\
g_2&=G\circ G_2\circ\mu_2\\
h_1&=\mu_1^{-1}\circ H_1\circ H\\
h_2&=\mu_2^{-1}\circ H_2\circ H
\end{align*}
for some $G\in\C[x]$, some $H\in\C(x)$, and some linear
$\mu_1,\mu_2\in\C[x]$, where $(G_1,G_2)$ satisfy one of
\emph{(\ref{BTcor}.1)}--\emph{(\ref{BTcor}.5)} and $(H_1,H_2)$ is the
corresponding pair below:
\begin{enumerate}
\item[(\thethm.1)]$(x^r p(x^n), \,x^n)$;
\item[(\thethm.2)]$((x-1/x)p(x+1/x), \,x+1/x)$;
\item[(\thethm.3)]$(D_n(x), \,D_m(x))$;
\item[(\thethm.4)]$\left(x^2+2x+\frac{1}{x}-\frac{1}{4x^2},
\,\frac{1}{3}\left((x+1-\frac{1}{2x})^3+4\right)\right)$;
\item[(\thethm.5)]$(x^n+1/x^n, \,(\zeta x)^m+1/(\zeta x)^m)$ where
$\zeta^{dmn}=-1$.
\end{enumerate}
\end{thm}

\begin{proof}
Since $g_1(x)-g_2(y)$ vanishes when $x=h_1(z)$ and $y=h_2(z)$, there is an
irreducible factor $E(x,y)$ of $g_1(x)-g_2(y)$ such that $E(h_1(z),h_2(z))=0$.
Here $E=0$ defines a genus-zero curve having at most two closed points lying
over $x=\infty$.  By Corollary~\ref{BTcor}, we have
$g_1=G\circ G_1\circ\mu_1$ and $g_2=G\circ G_2\circ\mu_2$ where
$G,\mu_1,\mu_2\in\C[x]$ with $\mu_i$ linear, and moreover
(perhaps after switching $g_1$ and $g_2$) we may choose $(G_1,G_2)$ to have
the form of one of (\ref{BTcor}.1)--(\ref{BTcor}.5).  Furthermore,
these choices can be made so that $E(x,y)$ divides
$G_1\circ\mu_1(x)-G_2\circ\mu_2(y)$.  As noted in Proposition~\ref{BTthm},
in cases (\ref{BTcor}.1)--(\ref{BTcor}.4) the polynomial $G_1(x)-G_2(y)$
is irreducible.
Thus, for any $H_1,H_2\in\C(x)$ satisfying
$G_1\circ H_1=G_2\circ H_2$ and
$\gcd(\deg(H_1),\deg(H_2))=1$, there exists $H\in\C(x)$ such that
$\mu_1\circ h_1 = H_1\circ H$ and $\mu_2\circ h_2 = H_2\circ H$.
Hence in these cases it suffices to exhibit one such pair
$(H_1,H_2)$, and visibly the pairs stated in the Theorem have the
required properties.

Henceforth suppose that $G_1=D_{dm}$ and $G_2=-D_{dn}$ with $m,n$ coprime
positive integers and $d>1$.
Let $G_1=a_1\circ b_1$ and $G_2=a_2\circ b_2$ be the decompositions
occurring in Proposition~\ref{Fried}.  Denoting by $\Omega$ the splitting
field of $a_1(x)-z$ over $\C(z)$, we see that $\deg(a_1)$ is the ramification
index in $\Omega/\C(z)$ of any place lying over $z=\infty$; but
(\ref{Fried}.2) implies the same description applies to $\deg(a_2)$, so
$a_1$ and $a_2$ have the same degree.  By Lemma~\ref{Dickson}, there
exist linear $\nu_1,\nu_2\in\C[x]$, and a divisor $e$ of $d$, such that
$a_1=D_e\circ\nu_1$ and $a_2=-D_e\circ\nu_2$ (and
$b_1=\nu_1^{-1}\circ D_{md/e}$ and $b_2=\nu_2^{-1}\circ D_{nd/e}$).
Since Proposition~\ref{Fried} holds for some linear $\nu_1,\nu_2$,
it follows that Proposition~\ref{Fried} holds for any arbitrarily chosen
linears $\nu_1,\nu_2$, so we may assume $\nu_1=\nu_2=x$.
A factorization of $a_1(x)-a_2(y)$ is given in Proposition~\ref{Dfac},
in terms of the polynomials
$A_{k,e}:=x^2 - xy\cdot 2\cos(\pi k/e) + y^2 - 4\sin^2(\pi k/e)$
where $1\le k<e$ and $k$ is odd.
Note that $A_{k,e}$ is irreducible (since its degree-$2$ part is a nonsquare,
it has no degree-$1$ terms, and it has a nonzero constant term).
Thus, by (\ref{Fried}.3), every irreducible factor of
$G_1(x)-G_2(y)$ has $x$-degree $2dm/e$, unless $e$ is odd when there is
also one factor of $x$-degree $dm/e$.  But Proposition~\ref{Dfac}
implies that $G_1(x)-G_2(y)$ is the product of several polynomials
$A_{k,d}(D_m(x),D_n(y))$, as well as (if $d$ is odd) the polynomial
$D_m(x)+D_n(y)$.  Thus every irreducible factor of $G_1(x)-G_2(y)$ has
$x$-degree at most $2m$, so either $e=d$ or $(d,e)=(2,1)$.
In the latter case, $G_1(x)-G_2(y)$ is irreducible.  Thus, in either case,
the irreducible factors of $G_1(x)-G_2(y)$ consist just of the
polynomials $A_{k,d}(D_m(x),D_n(y))$ with $1\le k<d$ and $k$ odd,
unless $d$ is odd in which case $D_m(x)+D_n(y)$ is another irreducible factor.
Now $E(\mu_1^{-1}(x),\mu_2^{-1}(y))$ must be a scalar multiple of one of
these factors, and we may assume the scalar is $1$ (since we are free to
replace $E$ by a scalar multiple of itself).
Since $E(h_1(x),h_2(y))=0$, we cannot have
$E(\mu_1^{-1}(x),\mu_2^{-1}(y))=D_m(x)+D_n(y)$, so we must have
$E(\mu_1^{-1}(x),\mu_2^{-1}(y))=A_{k,d}(D_m(x),D_n(y))$.
Denote this polynomial as $R(x,y)$, and
put $H_1:=x^n+1/x^n$ and $H_2:=(\zeta x)^m+1/(\zeta x)^m$,
where $\zeta=e^{\pi i k/(dmn)}$.  Then $R(H_1(x),H_2(x))=0$,
so (since $R(x,y)$ is irreducible) we have
$H_1=\hat{H_1}\circ J$ and $H_2=\hat{H_2}\circ J$ for some
$\hat{H_1},\hat{H_2},J\in\C(x)$ such that $R(\hat{H_1}(x),\hat{H_2}(x))=0$,
where in addition $\mu_1\circ h_1 = \hat{H_1}\circ H$ and
$\mu_2\circ h_2 = \hat{H_2}\circ H$ for some $H\in\C(x)$.
If $\deg(J)=1$ this gives (\ref{Type1}.5), so assume $\deg(J)>1$.
Since $\deg(J)$ divides $\gcd(\deg(H_1),\deg(H_2))=2$, we must have
$\deg(J)=2$.  If $J\in\C(x^2)$ then $H_1,H_2\in\C(x^2)$ so both $n$ and
$m$ are even, contradiction.  Now Lemma~\ref{Dickson} implies that
$J=\lambda_1\circ (x/\gamma+\gamma/x)$ and
$J=\lambda_2\circ (x/\delta+\delta/x)\circ\zeta x$, where
$\gamma^{2n}=1=\delta^{2m}$
and $\lambda_1,\lambda_2\in\C(x)$ have degree one.  Comparing images
of $x=0$, we see that $\lambda_1(\infty)=\lambda_2(\infty)$, so
$\lambda_2^{-1}\circ\lambda_1$ fixes $\infty$ and thus is a linear
polynomial.  Thus $J$ is a Laurent polynomial, and its constant term
is $\lambda_1(0)=\lambda_2(0)$, so $\lambda_2^{-1}\circ\lambda_1 = \epsilon x$
for some $\epsilon\in\C^*$.  Thus
\[
\epsilon\left(\frac{x}{\gamma}+\frac{\gamma}{x}\right) =
\frac{\zeta x}{\delta}+\frac{\delta}{\zeta x},
\]
and equating coefficients of like terms yields
$\epsilon\delta=\zeta\gamma$ and $\epsilon\gamma\zeta=\delta$, so
$\epsilon=\zeta\gamma/\delta=\pm 1$.  Raising to the $(2nm)\tth$ power gives
$\zeta^{2mn}=1$, but $\zeta^{2mn}=e^{2 \pi i k/d}\ne 1$ since $0< k<d$,
contradiction.
\end{proof}


\section{Laurent polynomials with decompositions of both types}
\label{sec-2}

In this section we prove the following result:

\begin{thm}\label{Type12}
Let $g_1\in\C[x]\setminus\C$ and $g_2,h_1\in\cL\setminus\C$ satisfy
$g_1\circ h_1=g_2\circ x^n$ with $n\in\bZ_{>0}$.  Then either $h_1=A\circ x^n$
(and $g_2=g_1\circ A$) for some $A\in\cL$, or there exist $G,\mu\in\C[x]$
with $\mu$ linear such that $g_1=G\circ G_1\circ\mu$ and
$h_1=\mu^{-1}\circ H_1$ and $g_2=G\circ G_2$, where one of the following holds
(with $e\in\bZ$ and $r=\gcd(n,e)$):
\begin{enumerate}
\item[(\thethm.1)] $G_1=x^{n/r}$, $\,H_1=x^e p(x^n)$, and
$G_2=x^{e/r} p(x)^{n/r}$, where $p\in\C[x]$;
\item[(\thethm.2)] $G_1=D_{n/r}$, $\,H_1=(x^e+1/x^e)\circ\alpha x$,
and $G_2=(x^{e/r}+1/x^{e/r})\circ\alpha^n x$, where $\alpha\in\C^*$.
\end{enumerate}
\end{thm}

We will use some results of Avanzi and Zannier \cite[\S 4]{AZ}, which we
state as follows.

\begin{prop}[Avanzi--Zannier]\label{AZprop}
Pick an indecomposable $g\in\C[x]$, and distinct nonconstant
$h_1,h_2\in\C(x)$, and suppose that $g\circ h_1=g\circ h_2$.
Then $g=\mu\circ G\circ\nu$ and $h_1=\nu^{-1}\circ H_1\circ H$ and
$h_2=\nu^{-1}\circ H_2\circ H$, where $\mu,\nu\in\C[x]$ are linear,
$H\in\C(x)$, and either $(G,H_1,H_2)$ or $(G,H_2,H_1)$ is in the
following list:
\begin{enumerate}
\item[(\thethm.1)] $(x^n, \,x, \,\zeta x)$, where
$n$ is prime and $\zeta$ is a primitive $n\tth$ root of unity;
\item[(\thethm.2)] $\left(D_n, \,x+\frac{1}{x},
 \,\zeta x+\frac{1}{\zeta x}\right)$, where $n$ is an odd prime
and $\zeta$ is a primitive $n\tth$ root of unity;
\item[(\thethm.3)] $\left((x^r(x+1)^m, \,\frac{1-x^r}{x^{r+m}-1},
\,-1+\frac{x^m-1}{x^{r+m}-1}\right)$, where $r,m$ are coprime positive
integers with $r+m>3$;
\item[(\thethm.4)] $\left(x(x+\alpha)^2(x+1)^2,
\,-4\alpha \frac{x^2}{E}, \,-\frac{\alpha}{E}\left(x^2-\frac{7x}{4}-
\frac{15}{64}\right)^2\right)$, where $\alpha\in\C^*$ satisfies
$9\alpha^2-2\alpha+9=0$ and
\[E=\alpha x^4 + \frac{3}{8}(3-7\alpha)x^3 +
 \frac{99}{64}(1+\alpha)x^2 + \frac{45}{512}(7-3\alpha)x +
    \frac{225}{4096};\]
\item[(\thethm.5)] $\left(x(x+\alpha)^3(x+1)^3, \,-4096\frac{x^3}{E},
\,\frac{1}{E}(64-(x-\alpha)^2)^3\right)$, where $\alpha\in\C^*$
satisfies $\alpha^2-5\alpha+8=0$ and
\begin{align*}
E=x^6 &+ (32-10\alpha)x^5 + (31\alpha - 88)x^4 + (68\alpha + 1888)x^3\\
&+ (651\alpha - 56)x^2 + (11158\alpha - 50288)x + 41881\alpha - 156520.
\end{align*}
\end{enumerate}
\end{prop}

\begin{remark}
The polynomials in \cite{AZ} involve some parameters
which we have removed by absorbing them into $\mu$ and $\nu$.
Also, the assertion in \cite[Prop.~4.7]{AZ} about $g_1$ being reduced
is false in case (3).
\end{remark}

\begin{proof}[Proof of Theorem~\ref{Type12}]
Let $\zeta$ be a primitive $n\tth$ root of unity, so for
$h_2:=h_1\circ \zeta x$ we have $g_1\circ h_2 = g_1\circ h_1$.
If $h_2=h_1$ then $h_1=A\circ x^n$ with $A\in\cL$, in which case
$g_2=g_1\circ A$.  Henceforth assume $h_2\ne h_1$.  This implies $g_1$ is
not linear, so we can write $g_1=f_1\circ\dots\circ f_v$ where every
$f_i$ is indecomposable.  Let $j$ be the largest integer for which
\[
f_j\circ f_{j+1}\circ\dots\circ f_v\circ h_2 =
f_j\circ f_{j+1}\circ\dots\circ f_v\circ h_1,
\]
and put $R=f_{j+1}\circ\dots\circ f_v$ and $A=f_1\circ\dots\circ f_{j-1}$,
so $g_1=A\circ f_j\circ R$.
Then $S_2:=R\circ h_2$ and $S_1:=R\circ h_1$ satisfy $S_2\ne S_1$ but
$f_j\circ S_2=f_j\circ S_1$.  After replacing $A$, $f_j$, and $R$ by
$A\circ\mu$, $\mu^{-1}\circ f_j\circ\nu^{-1}$, and $\nu\circ R$,
for suitable linear $\mu,\nu\in\bC[x]$, Proposition~\ref{AZprop} implies
that there exist $s_1,s_2,T\in\bC(x)\setminus\bC$ such that
$S_1=s_1\circ T$ and $S_2=s_2\circ T$ and either $(f_j,s_1,s_2)$ or
$(f_j,s_2,s_1)$ is one of the triples (\ref{AZprop}.1)--(\ref{AZprop}.5).
Since replacing $\zeta$ by $1/\zeta$ has the effect of exchanging $s_1$
and $s_2$, we may assume that $(f_j,s_1,s_2)$ is among
(\ref{AZprop}.1)--(\ref{AZprop}.5).
Moreover, since $h_i\in\cL$ and $R\in\C[x]$, also $S_i=R\circ h_i$ is in
$\cL$, so $s_i$ has at most two poles.  This rules out
(\ref{AZprop}.3), (\ref{AZprop}.4) and (\ref{AZprop}.5).

If (\ref{AZprop}.1) holds then $f_j=x^{\ell}$ for some prime $\ell$,
and moreover $S_2=\gamma S_1$ for some primitive $\ell\tth$ root of unity
$\gamma$.  Thus $S_1(\zeta x) = \gamma S_1(x)$, so $S_1\in x^t\C[x^n]$ for some
$t\in\bZ$ with $\zeta^t=\gamma$.  By Proposition~\ref{Ritttwisthard},
after replacing $R$ and $h_1$ by $R\circ\mu$ and $\mu^{-1}\circ h_1$ for
a suitable linear $\mu\in\C[x]$, we may assume that $R$ and $h_1$ satisfy
the conditions required of $g$ and $h$ in either (\ref{Ritttwisthard}.1) or
(\ref{Ritttwisthard}.2).
First suppose $R$ and $h_1$ satisfy (\ref{Ritttwisthard}.1), so
$R\in x^d\C[x^m]$ and $h_1=x^e p(x^n)$ with $p\in\C[x]$ and
$n\mid em$; then $\gamma=\zeta^{de}$, so $n\mid \ell de$.
Putting $r=\gcd(n,e)$, we have $n\mid r\gcd(d\ell,m)$, so
$f_j\circ R\in\C[x^{n/r}]$.  Since $g_1=A\circ f_j\circ R$, we can write
$g_1=G\circ x^{n/r}$.  It follows that $g_2=x^{e/r} p(x)^{n/r}$, so we
have (\ref{Type12}.1).
Now suppose $R$ and $h_1$ satisfy (\ref{Ritttwisthard}.2), so $n$ is even,
and also $R=\hat{R}\circ D_d$ and
$h_1=(x^e+1/x^e)\circ\alpha x$, where $\hat{R}\in x\C[x^2]$ and
$ed\equiv t\equiv n/2\pmod{n}$; thus $\gamma=\zeta^t$ has
order $2$, so $\ell=2$.
Now $f_j\circ R = x^2\circ\hat{R}\circ D_d$; since $\hat{R}\in x\C[x^2]$,
we see that $x^2\circ\hat{R}$ is in $\C[x^2]$, and thus can be written as
$\widetilde{R}\circ D_2$ with $\widetilde{R}\in\C[x]$.
Thus $f_j\circ R=\widetilde{R}\circ D_{2d}$,
so since $n\mid 2ed$ we can write $g_1=G\circ D_{n/r}$ where $r=\gcd(n,e)$
amd $G=A\circ\widetilde{R}\circ D_{2dr/n}$.  This implies
$g_2=(x^{e/r}+x^{-e/r})\circ\alpha^n x$, so we have (\ref{Type12}.2).

Finally, suppose (\ref{AZprop}.2) holds.  Then $f_j=D_{\ell}$ for some
odd prime $\ell$, and moreover
$s_1=x+1/x$ and $s_2=\gamma x+1/(\gamma x)$ for some primitive $\ell\tth$
root of unity $\gamma$.  Since $s_1\circ T$ is a Laurent polynomial,
and $s_1$ has poles at $0$ and $\infty$, Lemma~\ref{basic} implies
that $T=\delta x^d$ for
some $\delta\in\bC^*$ and $d\in\bZ$.  Since we can replace $s_1$, $s_2$,
and $T$ by $s_1\circ 1/x$, $s_2\circ 1/x$, and $1/x\circ T$, we may
assume $d>0$.  Now we have $R\circ h_1=\delta x^d+1/(\delta x^d)$ and
$R\circ h_1(\zeta x) = \gamma \delta x^d+1/(\gamma \delta x^d)$, so
$\zeta^d=\gamma$ and thus $n\mid d\ell$.
Since $R$ is a polynomial, Lemma~\ref{Dickson} implies that
$R=\alpha^d D_{d/e}\circ\mu$ and
$h_1=\mu^{-1}\circ (x^e+1/x^e)\circ \hat\delta x/\alpha$
where $\mu\in\C[x]$ is linear, $\alpha^{2d}=1$, and $\hat\delta^d=\delta$.
Likewise
$R=\beta^d D_{d/e}\circ\nu$ and $h_1\circ\zeta x = \nu^{-1}\circ
(x^e+1/x^e)\circ \hat\gamma x/\beta$ where $\nu\in\C[x]$ is linear,
$\beta^{2d}=1$, and $\hat\gamma^d=\gamma \delta$.
Thus $(\alpha/\beta)^d D_{d/e} = D_{d/e}\circ
\nu\circ\mu^{-1}$; equating coefficients of $x^{d/e-1}$ shows that
$\nu\circ\mu^{-1}=\theta x$ with $\theta\in\C^*$, and equating coefficients of
$x^{d/e}$ shows that $\theta^{d/e}=(\alpha/\beta)^d$.  If $d=e$ it follows that
$\theta\in\{1,-1\}$; if $d\ne e$ then we also obtain $\theta=\pm 1$ upon
equating coefficients of $x^{d/e-2}$.
Since $\epsilon:=\alpha^d=\pm 1$, we have
\begin{align*}
g_1&=A\circ D_{\ell}\circ \epsilon D_{d/e}\circ\mu \\
&=A\circ \epsilon^{\ell} D_{\ell}\circ D_{d/e}\circ\mu
 \quad\text{(by (\ref{cD}))}\\
&=A\circ\epsilon^{\ell} D_{\ell d/e}\circ\mu.
\end{align*}
Recall that $n\mid d\ell$, so with $r=\gcd(e,n)$ we have
$en\mid \ell dr$, and thus
$g_1=G\circ D_{n/r}\circ\mu$ with $G=A\circ\epsilon^{\ell} D_{\ell d r/(en)}$.
Since $h_1=\mu^{-1}\circ (x^e+1/x^e)\circ \hat\delta x/\alpha$,
we find $g_2=G\circ (x^{e/r}+x^{-e/r})\circ (\hat\delta/\alpha)^nx$,
so we have (\ref{Type12}.2).
\end{proof}


\section{Proofs of main results}
\label{sec-proofs}

In this section we prove the results stated in Section~1.


\subsection{Proof of Theorem~\ref{LRitt1}}
Define an `admissible sequence' to be a finite sequence of complete
decompositions of a rational function $f$, such that consecutive decompositions
in the sequence differ only in that two adjacent indecomposables $u,v$ in the
first decomposition are replaced in the second decomposition by two other
indecomposables $\hat{u},\hat{v}$ such that $u\circ v=\hat{u}\circ\hat{v}$
and $\{\deg(u),\deg(v)\}=\{\deg(\hat{u}),\deg(\hat{v})\}$.
It suffices to prove that, for any two complete decompositions of a Laurent
polynomial $f$, there is an admissible sequence containing them both.
We prove this by induction on $\deg(f)$.  So assume it holds for all
Laurent polynomials of degree less than $\deg(f)$, and consider two
complete decompositions
$f=p_1\circ p_2\circ\dots\circ p_r=q_1\circ q_2\circ\dots\circ q_s$
(so $p_i,q_j\in\C(x)$ are indecomposable).  If $r=1$ or $s=1$ then these
decompositions are identical, so trivially are contained in an
admissible sequence.  Henceforth assume $r,s>1$.

By Lemma~\ref{basic}, after replacing $p_{r-1}$ and $p_r$ by
$p_{r-1}\circ\mu$ and $\mu^{-1}\circ p_r$ for some $\mu\in\C(x)$ with
$\deg(\mu)=1$, we may assume that
both $p_r$ and $\hat{p}:=p_1\circ\dots\circ p_{r-1}$ are Laurent polynomials,
and moreover either $\hat{p}\in\C[x]$ or $p_r=x^n$ with $n$ prime.
Further, if $\hat{p}\in\C[x]$ and $p_r\in\C(x^n)$ with $n>1$, then
$n$ is prime and $p_r=\hat\mu\circ x^n$ for some degree-one
$\hat\mu\in\C(x)$, so by replacing $p_{r-1}$ and $p_r$ by
$p_{r-1}\circ\hat\mu$ and $\hat\mu^{-1}\circ p_r$ we may assume
$p_r=x^n$; since $\hat{p}\circ p_r=f\in\cL$, we must have $\hat{p}\in\cL$.
Thus we may assume that $\hat{p},p_r\in\cL$,
and if there is no prime $n$ for which $p_r=x^n$, then $\hat{p}\in\C[x]$ and
$p_r\notin\C(x^n)$ for any $n>1$.  We can make analogous assumptions about
$q_s$ and $\hat{q}:=q_1\circ\dots\circ q_{s-1}$.

If there is a degree-one $\nu\in\C(x)$ for which $p_r=\nu\circ q_s$, then
$\hat{p}=\hat{q}\circ\nu^{-1}$, so by induction there is an admissible
sequence containing $p_1\circ\dots\circ p_{r-1}$ and
$q_1\circ\dots\circ q_{s-2}\circ (q_{s-1}\circ\nu^{-1})$.  Composing
each complete decomposition in the sequence with $p_r$, we then get an
admissible sequence containing $p_1\circ\dots\circ p_r$ and
$q_1\circ\dots\circ q_s$.  Henceforth assume there is no such $\nu$.

If $p_r=x^n$ and $q_s=x^m$ (with $n,m$ distinct primes), then
Proposition~\ref{Type2} implies
$\hat{p}=G\circ x^m$ and
$\hat{q}=G\circ x^n$ for some $G\in\cL$. Write
$G=g_1\circ\dots\circ g_t$ where every
$g_i\in\C(x)$ is indecomposable.
By induction, there is an admissible sequence containing
$p_1\circ\dots\circ p_{r-1}$ and
$g_1\circ\dots\circ g_t\circ x^m$, so composing with $p_r$ yields an
admissible sequence containing
$p_1\circ\dots\circ p_{r-1}\circ p_r$
and $g_1\circ\dots\circ g_t\circ x^m\circ x^n$.
Likewise there is an admissible sequence containing
$q_1\circ\dots\circ q_s$ and
$g_1\circ\dots\circ g_t\circ x^n\circ x^m$.
Since the sequence $(x^m\circ x^n, x^n\circ x^m)$ is admissible,
there is an admissible sequence containing
$p_1\circ\dots\circ p_r$ and $q_1\circ\dots\circ q_s$.

Now assume $q_s=x^n$ but $p_r\notin\C(x^m)$ for every $m>1$.
Then $\hat{p}\in\C[x]$.  By Theorem~\ref{Type12},
there exist $G,\mu\in\C[x]$ with $\deg(\mu)=1$
such that $\hat{p}=G\circ G_1\circ\mu$ and
$p_r=\mu^{-1}\circ H_1$ and $\hat{q}=G\circ G_2$, where
$G_1,G_2,H_1$ satisfy either (\ref{Type12}.1) or (\ref{Type12}.2).
In (\ref{Type12}.2) we have $H_1 = (x^e+1/x^e)\circ\alpha x$ with
$\alpha\in\C^*$ and $e>0$, and indecomposability of $p_r$
implies $e=1$.  Thus $G_1=D_n$ and $G_2=(x+1/x)\circ\alpha^n x$,
so $(G_1\circ H_1, G_2\circ q_s)$ is admissible, and the inductive
argument of the previous paragraph produces an admissible sequence
containing $p_1\circ\dots\circ p_r$ and $q_1\circ\dots\circ q_s$.
In (\ref{Type12}.1) we have $H_1=x^e h(x^n)$ with $h\in\C[x]$ and
$e\in\bZ$; since $p_r\notin\C(x^m)$ for $m>1$, we must have $\gcd(e,n)=1$,
so $G_1=x^n$ and $G_2=x^e h(x)^n$.  We will show that $G_2$ is indecomposable.
This implies that $(G_1\circ H_1, G_2\circ x^n)$ is admissible, so as above
there is an admissible sequence containing $p_1\circ\dots\circ p_r$ and
$q_1\circ\dots\circ q_s$.  So suppose $G_2$ is decomposable; then
Lemma~\ref{basic} implies that $G_2$ has a decomposition of either
Type 1 or Type 2 in which both rational functions involved have
degree $>1$.  By Proposition~\ref{Ritttwisteasy}, if there is a Type 1
decomposition with this property, then $G_2=u\circ v$ where
$u=x^iA(x)^n$
and $v=x^jB(x)^n$, with $A,B\in\C[x]$ and $i,j\in\bZ$ and $i>0$.
But then $x^n\circ H_1 = G_2(x^n) = u\circ v(x^n) =
x^iA^n\circ x^n\circ x^jB(x^n) =
x^n\circ x^i A(x^n)\circ x^j B(x^n)$, so
$H_1=\zeta x^i A(x^n)\circ x^j B(x^n)$ for some $\zeta\in\C^*$
with $\zeta^n=1$, contradicting indecomposability of $p_r$.
If $G_2$ has a Type 2 decomposition into rational functions of degree $>1$,
say $G_2\in\C(x^m)$ with $m>1$, then $G_2(\zeta x)=G_2(x)$
where $\zeta$ is a primitive $m\tth$ root of unity.  Thus
$\zeta^e h(\zeta x)^n = h(x)^n$, so $h(\zeta x) = \beta h(x)$ where
$\zeta^e \beta^n = 1$.  Hence $h=x^d A(x^m)$ for some $A\in\C[x]$
and some $d\in\bZ$ such that $\zeta^d=\beta$.  Thus
$1=\zeta^e \beta^n = \zeta^{e+nd}$, so $m\mid (e+nd)$.  Now
$H_1 = x^e h(x^n) = x^{e+nd} A(x^{nm})$ is in $\C(x^m)$, and since
$H_1$ is indecomposable we must have $H_1=\lambda\circ x^m$ for some
degree-one $\lambda\in\C(x)$.  But $H_1=x^e h(x^n)$ has no constant
term (since $\gcd(e,n)=1$), so $\lambda$ is a degree-one Laurent polynomial
with no constant term, whence $\lambda$ is a monomial Laurent polynomial.
Thus $h$ is a monomial polynomial, so $G_2=x^e h(x)^n$ is a constant
times $H_1=x^e h(x^n)$, whence indecomposability of $H_1$ implies
indecomposability of $G_2$.

Now assume $p_r,q_s\notin\C(x^n)$ for every $n>1$.
This implies $\hat{p}, \hat{q}\in\C[x]$, so Theorem~\ref{Type1}
applies.  After switching $(\hat{p},p_r)$ and $(\hat{q},q_s)$ if
necessary, we obtain
\begin{align*}
\hat{p}&=G\circ G_1\circ\mu_1\\
\hat{q}&=G\circ G_2\circ\mu_2\\
p_r&=\mu_1^{-1}\circ H_1\circ H\\
q_s&=\mu_2^{-1}\circ H_2\circ H
\end{align*}
for some $H\in\C(x)$ and $G,\mu_1,\mu_2\in\C[x]$ with $\mu_i$ linear,
where $(G_1,G_2)$ is one of (\ref{BTcor}.1)--(\ref{BTcor}.5) and
$(H_1,H_2)$ is the corresponding pair among (\ref{Type1}.1)--(\ref{Type1}.5).
If $\deg(H)>1$ then indecomposability of $p_r$ and $q_s$ implies
$p_r=\nu\circ q_s$ for some degree-one $\nu\in\C(x)$, a case treated
previously.  So assume $\deg(H)=1$, whence $H_1$ and $H_2$ are
indecomposable.  In case (\ref{Type1}.1) we have $H_2=x^n$ with $n>0$ (where
indecomposability implies $n$ is prime), and $H_1=x^e h(x^n)$ with
$h\in\C[x]$ and $e\in\bZ_{>0}$
coprime to $n$.  Moreover, $G_1=x^n$ and $G_2=x^e h(x)^n$.
Here indecomposability of $H_1$ implies indecomposability of $G_2$
(by Ritt's first theorem), so our result follows by induction.
In case (\ref{Type1}.2) we have $H_2=x+1/x$ and $H_1=(x-1/x)p(x+1/x)$
with $p\in\C[x]$, and moreover $G_1=x^2$ and $G_2=(x^2-4)p(x)^2$.
Here we need only to prove that $G_2$ is indecomposable.  If it were
not, then by Proposition~\ref{decompweird1} there would be nonlinear
$u,v\in\C[x]$ such that $u\circ v = G_2$ and $u,v$ satisfy the conditions
required of $g,h$ in either (\ref{decompweird1}.1)
or (\ref{decompweird1}.2).  In (\ref{decompweird1}.1) we have
$u=xB^2$ and $v=(x^2-4)D^2$ with $B,D\in\C[x]$, so composing with $x+1/x$
gives
\begin{align*}
x^2\circ H_1 &= G_2\left(x+\frac{1}{x}\right) = u\circ v\left(x+\frac{1}{x}\right)
= u\circ x^2\circ \left(x-\frac{1}{x}\right)\cdot D\left(x+\frac{1}{x}\right) \\
&= x^2\circ xB(x^2)\circ \left(x-\frac{1}{x}\right)\cdot D\left(x+\frac{1}{x}\right),
\end{align*}
whence $H_1 = \pm xB(x^2)\circ (x-1/x)D(x+1/x)$, contradicting indecomposability
of $H_1$.  In (\ref{decompweird1}.2) we have
$u=(x^2-4)B^2$ and $v=D_n$ where $B\in\C[x]$ and $n>1$, so composing with $x+1/x$
gives 
\begin{align*}
x^2\circ H_1 &= G_2\left(x+\frac{1}{x}\right) = u\circ v\left(x+\frac{1}{x}\right) =
 u\circ \left(x+\frac{1}{x}\right)\circ x^n \\
&= x^2\circ \left(x-\frac{1}{x}\right)\cdot B\left(x+\frac{1}{x}\right)\circ x^n,
\end{align*}
whence $H_1=\pm (x-1/x)B(x+1/x)\circ x^n$,
again contradicting indecomposability.  If (\ref{Type1}.3) holds then $H_1=G_2=D_n$
and $H_2=G_1=D_m$ where $m,n$ are distinct primes, so the result follows by
induction.
If (\ref{Type1}.4) holds then $H_2$ is decomposable, a contradiction.
Suppose (\ref{Type1}.5) holds.  Then $H_1=x^n+1/x^n$ with $n\in\bZ_{>0}$, and
indecomposability implies $n=1$.  Likewise $H_2=\zeta x+1/(\zeta x)$, where
$\zeta^d=-1$ for some $d\in\bZ_{>1}$, and moreover $G_1=D_d=-G_2$.  Write
$d=\prod_{i=1}^t \ell_i$ where the $\ell_i$ are primes which need not be
distinct, and put $e=d/\ell_1$.  Since
$D_e\circ (x+1/x) = -D_e\circ (\zeta^{\ell_1} x+1/(\zeta^{\ell_1} x))$, by induction
there is an admissible sequence containing both
$D_{\ell_2}\circ\dots\circ D_{\ell_t}\circ (x+1/x)$ and
$-D_{\ell_2}\circ D_{\ell_3}\circ\dots\circ D_{\ell_t}\circ
 (\zeta^{\ell_1} x+1/(\zeta^{\ell_1} x))$.  Composing with $x^{\ell_1}$ gives an
admissible sequence containing
$D_{\ell_t}\circ\dots\circ D_{\ell_2}\circ (x+1/x)\circ x^{\ell_1}$ and
$-D_{\ell_t}\circ D_{\ell_{t-1}}\circ\dots\circ D_{\ell_2}\circ
 (\zeta^{\ell_1} x+1/(\zeta^{\ell_1} x))\circ x^{\ell_1}$, and plainly
$((x+1/x)\circ x^{\ell_1}, \,D_{\ell_1}\circ (x+1/x))$ is admissible, as is
$((\zeta^{\ell_1}x+1/(\zeta^{\ell_1}x))\circ x^{\ell_1}, \,
 D_{\ell_1}\circ (\zeta x+1/(\zeta x)))$.
Thus there is an admissible sequence containing
$D_{\ell_t}\circ\dots\circ D_{\ell_1}\circ H_1$ and
$-D_{\ell_t}\circ D_{\ell_{t-1}}\circ\dots\circ D_1\circ H_2$,
so there is an admissible sequence containing
$p_1\circ\dots\circ p_r$ and $q_1\circ\dots\circ q_s$.
This concludes the proof of Theorem~\ref{LRitt1}.


\subsection{Proof of Theorem~\ref{LRitt2}}
We prove the result by induction on $\deg(f)$.  So assume it holds
for all Laurent polynomials of degree less than $\deg(f)$, and write
$f=g_1\circ h_1 = g_2\circ h_2$ with $f\in\cL$ and with indecomposable
$g_1,g_2,h_1,h_2\in\C(x)$.
After replacing $g_1$ and $h_1$ by $g_1\circ\mu$ and
$\mu^{-1}\circ h_1$ for some $\mu\in\C(x)$ with $\deg(\mu)=1$, we may
assume that $g_1,h_1\in\cL$ and either $g_1\in\C[x]$ or $h_1=x^n$
with $n$ prime (by Lemma~\ref{basic}).
Moreover, this argument shows that if $h_1\in\C(x^n)$ for some $n>1$
we may assume $h_1=x^n$ (and indecomposability implies $n$ is prime).
We can make analogous assumptions about $g_2$ and $h_2$.
If $h_1=\mu\circ h_2$ for some degree-one $\mu\in\C(x)$, then
$g_1\circ\mu=g_2$, so we have (\ref{LRitt2}.1).  Henceforth assume
$h_1\ne \mu\circ h_2$ for any degree-$1$ $\mu\in\C(x)$.

First suppose $h_1=x^m$ and $h_2=x^n$, where $m$ and $n$ are distinct
primes.  Proposition~\ref{Type2} implies
$g_1=G\circ x^n$ and
$g_2=G\circ x^m$ for some $G\in\cL$, which must have degree $1$
since $g_1$ and $h_2$ are indecomposable.  This yields (\ref{LRitt2}.2) with
$r=m$ and $q=1$.

Now suppose $h_2=x^n$ but $h_1\notin\C(x^m)$ for any $m>1$.
Then $g_1\in\C[x]$, so Theorem~\ref{Type12}
applies.  Since $h_1\notin\C(x^n)$, there
exist $G,\mu\in\C[x]$ with $\deg(\mu)=1$ such that 
$g_1=G\circ G_1\circ\mu$ and
$h_1=\mu^{-1}\circ H_1$ and $g_2=G\circ G_2$, where
either (\ref{Type12}.1) or (\ref{Type12}.2) holds.
If $\deg(G)>1$ then indecomposability of $g_1$ implies $\deg(G_1)=1$, so
in both (\ref{Type12}.1) and (\ref{Type12}.2) we have $n\mid e$ and thus
$G_1\in\C(x^n)$, contradiction.
Hence $\deg(G)=1$, so $G_1$ is indecomposable and thus
$\gcd(n,e)=1$.
In (\ref{Type12}.1) we have $G_1=x^n$ and $H_1=x^e q(x^n)$ and
$H_2=x^e q(x)^n$, with $q\in\C[x]$ and $e\in\bZ$ coprime to $n$;
this gives (\ref{LRitt2}.2).
In (\ref{Type12}.2) we have $G_1=D_n$ and $H_1=(\alpha x)^e + 1/(\alpha x)^e$
and $G_2=(\alpha^n x)^e + 1/(\alpha^n x)^e$ with $\alpha\in\C^*$ and
$e\in\bZ$, and indecomposability implies $e=\pm 1$.  After adjusting
$G_1,G_2,H_1,H_2$ by composing with linears, this gives (\ref{LRitt2}.4).

Henceforth assume $h_1,h_2\notin\C(x^m)$ for every $m>1$.
Then $g_1,g_2\in\C[x]$, so Theorem~\ref{Type1} applies.  Thus,
after switching $(g_1,h_1)$ and $(g_2,h_2)$ if necessary, we have
\begin{align*}
g_1&=G\circ G_1\circ\mu_1\\
g_2&=G\circ G_2\circ\mu_2\\
h_1&=\mu_1^{-1}\circ H_1\circ H\\
h_2&=\mu_2^{-1}\circ H_2\circ H
\end{align*}
for some $G\in\C[x]$, some linear $\mu_1,\mu_2\in\C[x]$, and some
$H\in\C(x)$, where $(G_1,G_2)$ is one of the pairs
(\ref{BTcor}.1)--(\ref{BTcor}.5) and $(H_1,H_2)$ is the corresponding pair
among (\ref{Type1}.1)--(\ref{Type1}.5).
If $\deg(G)>1$ then indecomposability of $g_i$ implies $\deg(G_i)=1$,
so we must have either (\ref{BTcor}.1) or (\ref{Type1}.3).  Thus $(H_1,H_2)$
satisfy (\ref{Type1}.1) and (\ref{Type1}.3), and in either case
$G_1\circ H_1 = G_2\circ H_2$ is a linear polynomial, so we have
(\ref{LRitt2}.1).  Likewise if $\deg(H)>1$ then $\deg(H_i)=1$, so since
$g_1\circ\mu_1^{-1}\circ H_1 = g_2\circ\mu_2^{-1}\circ H_2$ we again
have (\ref{LRitt2}.1).  Now assume $\deg(G)=\deg(H)=1$, so
$G_i$ and $H_i$ are indecomposable.  Since $H_2\ne x^n$, we do not
have (\ref{Type1}.1).  If (\ref{BTcor}.2) and (\ref{Type1}.2) hold then, by
(\ref{q})--(\ref{q3}), there are $\nu_1,\nu_2,q\in\C(x)$ with
$\deg(\nu_i)=1$ such that $H_1\circ\nu_1=x q(x^2)$ and
$G_2\circ\nu_2=x q(x)^2$; here also $G_1=x^2$ and
$\nu_2^{-1}\circ H_2\circ \nu_1 = x^2$, so we have (\ref{LRitt2}.2).
Note that in this case $q$ is not a Laurent polynomial, instead $q=Q(1/(x+1))$
for some $Q\in x\C[x]$.
If (\ref{BTcor}.3) and (\ref{Type1}.3) hold then (\ref{LRitt2}.3) holds.
Since $G_1$ is indecomposable, we do not have (\ref{BTcor}.4).
Now suppose (\ref{BTcor}.5) and (\ref{Type1}.5) hold.
Thus $G_1=D_{dm}$ and $G_2=-D_{dn}$ with $d>1$ and $m,n\ge 1$, so
indecomposability implies $d$ is prime and $m=n=1$.  Here $H_1=x+1/x$ and
$H_2=H_1\circ\zeta x$, where $\zeta^d=-1$.  If $d$ is odd then,
with $\mu=-x$, we have $G_2=D_{dn}\circ\mu$ and
$\mu^{-1}\circ H_2=H_1\circ (-\zeta x)$ where $(-\zeta)^d=1$,
which is (\ref{LRitt2}.5).  Finally, if $d=2$ then with $\mu=2-x$ we see
that $(\mu\circ G_2,\,\mu\circ G_1)$ satisfies (\ref{BTcor}.2) and
$(H_2,H_1)$ satisfies (\ref{Type1}.2) (both with $p(x)=\zeta$),
a case we have already resolved.  This concludes the proof of
Theorem~\ref{LRitt2}.


\subsection{Proof of Theorem~\ref{Lbidec}}
Let $f\in\cL\setminus\C$ and $g_1,g_2,h_1,h_2\in\bC(x)$ satisfy
$f=g_1\circ h_1 = g_2\circ h_2$.  By Lemma~\ref{basic}, after replacing
$g_1$ and $h_1$ by $g_1\circ\mu$ and $\mu^{-1}\circ h_1$ for some
degree-one $\mu\in\C(x)$, we may assume $g_1,h_1\in\cL$ and either
$g_1\in\C[x]$ or $h_1=x^n$ with $n\in\bZ_{>0}$.  We can make similar
assumptions about $g_2$ and $h_2$.

If $h_1=x^n$ and $h_2=x^m$ with $n,m>0$, then Proposition~\ref{Type2}
implies $g_1=G\circ x^{\lcm(n,m)/n}$ and $g_2=G\circ x^{\lcm(n,m)/m}$
for some $G\in\cL$.  Thus (\ref{Lbidec}.1) holds with $\mu_i=x$ and
$H=x^{\gcd(n,m)}$ (and $p=1$).

Now suppose precisely one of $h_1$ and $h_2$ has the form $x^n$ with
$n>0$; by switching $(g_1,h_1)$ and $(g_2,h_2)$ if necessary, we
may assume $h_2=x^n$ and $g_1\in\C[x]$.
If there exists $A\in\cL$ such that
$h_1=A\circ x^n$ and $g_2=g_1\circ A$, then (\ref{Lbidec}.1) holds with
$G=g_1$, $\,\mu_i=x$, $\,H=x^n$, and $p=A$.
So assume there is no such $A$.  By Theorem~\ref{Type12}, there exist
$G,\mu\in\C[x]$
with $\mu$ linear such that $g_1=G\circ G_1\circ\mu$ and
$h_1=\mu^{-1}\circ H_1$ and $g_2=G\circ G_2$, where either
(\ref{Type12}.1) or (\ref{Type12}.2) holds.  If (\ref{Type12}.1) holds
then (\ref{Lbidec}.1) holds with $H=x^{\gcd(n,e)}$.
If (\ref{Type12}.2) holds then (\ref{Lbidec}.6) holds with
$H=(\alpha x)^{\gcd(n,e)}$.

Finally, suppose $g_1,g_2\in\C[x]$, so Theorem~\ref{Type1} applies.
Thus, perhaps after switching $(g_1,g_2)$ and $(h_1,h_2)$, we have
\begin{align*}
g_1&=G\circ G_1\circ\mu_1\\
g_2&=G\circ G_2\circ\mu_2\\
h_1&=\mu_1^{-1}\circ H_1\circ H\\
h_2&=\mu_2^{-1}\circ H_2\circ H
\end{align*}
for some $G\in\C[x]$, some $H\in\C(x)$, and some linear
$\mu_1,\mu_2\in\C[x]$, where $(G_1,G_2)$ satisfy one of
(\ref{BTcor}.1)--(\ref{BTcor}.5) and $(H_1,H_2)$ is the
corresponding pair among (\ref{Type1}.1)--(\ref{Type1}.5).
In each case, this implies the corresponding condition among
(\ref{Lbidec}.1)--(\ref{Lbidec}.5).

If $G_1\circ H_1$ has poles at both $0$ and
$\infty$, then $H$ preserves $\{0,\infty\}$, so $H$ is a monomial.
This occurs in (\ref{Lbidec}.2) and (\ref{Lbidec}.4)--(\ref{Lbidec}.6).

Now we prove the final assertion in Theorem~\ref{Lbidec}.
Since $f=G\circ G_1\circ H_1\circ H$ is a nonconstant Laurent polynomial,
and $G,G_1\in\C[x]$, we see that $H_1\circ H$ has no poles besides
$0$ and $\infty$.  If any of (\ref{Lbidec}.2) or
(\ref{Lbidec}.4)--(\ref{Lbidec}.6) holds, then $H_1$ has poles at both
$0$ and $\infty$, so $H$ preserves $\{0,\infty\}$ and thus $H=\alpha x^s$
with $\alpha\in\C^*$ and $s\in\bZ$.  Here $s\ne 0$ (since $f$ is nonconstant).
To show we can choose $s>0$, it suffices to prove that, for some
$\beta\in\C^*$ and some degree-one $\nu_1,\nu_2\in\C(x)$, the decompositions
$(G_1\circ\nu_1)\circ (\nu_1^{-1}\circ H_1\circ \beta/x) =
 (G_2\circ\nu_2)\circ (\nu_2^{-1}\circ H_2\circ \beta/x)$
satisfy the same one of (\ref{Lbidec}.2) or (\ref{Lbidec}.4)--(\ref{Lbidec}.6)
that is satisfies by the original decompositions.
In case (\ref{Lbidec}.2) this is true for $\beta=1$ and
$\nu_2=x=-\nu_1$.  
In (\ref{Lbidec}.4), we can take $\beta=-1/2$ and $\nu_2=x=-\nu_1$.
In (\ref{Lbidec}.5), we can take $\beta=1$ and $\nu_1=x=\nu_2$
(provided we replace $\zeta$ by $1/\zeta$).
In (\ref{Lbidec}.6), we can take $\beta=1$ and $\nu_1=x=1/\nu_2$.
This concludes the proof of Theorem~\ref{Lbidec}.


\subsection{Proof of Proposition~\ref{Lunique}}
Pick $f\in\cL\setminus\C$, and suppose there are $g_1,g_2,h_1,h_2\in\C(x)$
such that $f=g_1\circ h_1 = g_2\circ h_2$ and $\deg(g_1)=\deg(g_2)$.
By Theorem~\ref{Lbidec}, after possibly switching $(g_1,h_1)$ and
$(g_2,h_2)$, we have
\begin{align*}
g_1&=G\circ G_1\circ\mu_1\\
g_2&=G\circ G_2\circ\mu_2\\
h_1&=\mu_1^{-1}\circ H_1\circ H\\
h_2&=\mu_2^{-1}\circ H_2\circ H
\end{align*}
for some $G\in\C[x]$, some $H\in\cL$, and some degree-one
$\mu_1,\mu_2\in\C(x)$, where one of (\ref{Lbidec}.1)--(\ref{Lbidec}.6) holds.
Since $\deg(g_1)=\deg(g_2)$, we have $\deg(G_1)=\deg(G_2)$, which greatly
restricts the possibilities.  In particular, (\ref{Lbidec}.4) cannot happen.
In case (\ref{Lbidec}.3) we must have $m=n=1$, so (\ref{Lunique}.1) holds.
In case (\ref{Lbidec}.5) we again have $m=n=1$, so (\ref{Lunique}.3) holds.
In case (\ref{Lbidec}.6) we have $m=2$ and $n=1$, so (\ref{Lunique}.4) holds.
In case (\ref{Lbidec}.2) we have $p=\alpha x$ with $\alpha\in\C^*$.  Putting
$\lambda=2+\frac{x}{\alpha^2}$ and $\nu=\alpha x/i$ we get
\begin{align*}
\lambda\circ G_1\circ \nu &= -D_2(x) \\
\lambda\circ G_2 &= D_2(x) \\
\nu^{-1}\circ H_1 &= ix+\frac{1}{ix} \\
H_2 &= x+\frac{1}{x},
\end{align*}
which is the $n=2$ case of (\ref{Lunique}.3).
Finally, suppose (\ref{Lbidec}.1) holds, so
$G_1=H_2=x^n$ for some $n>0$, and $H_1=x^r p(x^n)$ 
and $G_2=x^r p(x)^n$ where $p\in\C[x]\setminus\{0\}$ and
$r\in\bZ$ is coprime to $n$.  Write $p=x^e P$ where $P\in\C[x]$
satisfies $P(0)\ne 0$, so with $R=r-en$ we have
$H_1=x^R P(x^n)$ and $G_2=x^R P(x)^n$; replacing $r$ by $R$ and
$p$ by $P$, we may therefore assume $x\nmid p$.  
If $r\ge 0$ then $\deg(G_2)=r+n\cdot\deg(p)$, which must equal $n$,
so $\deg(p)\le 1$.  In either case, coprimality of $r$ and $n$ implies
$n=1$: for, if $\deg(p)=1$ then $r=0$, and if $\deg(p)=0$ then $r=n$.
Thus $G_2$ and $H_1$ are linear, and $G_2=H_1$, so by composing with
linears we obtain (\ref{Lunique}.1).  Now assume $r<0$, and write $s=-r$.
Then $\deg(G_2)=\max(s,n\deg(p))$, so $\deg(p)=1$ and $1\le s\le n$.
We may assume $s<n$, since otherwise $s=n=1$ so we obtain
(\ref{Lunique}.1) as above.
Now, composing with (scalar) linears gives (\ref{Lunique}.2).

Cases (\ref{Lunique}.3) and (\ref{Lunique}.4) are instances of
(\ref{Lbidec}.5) and (\ref{Lbidec}.6), so by Theorem~\ref{Lbidec} we may
assume $H=\alpha x^s$ with $\alpha\in\C^*$ and $s\in\bZ_{>0}$.
If (\ref{Lunique}.2) holds, then $f=G\circ G_1\circ H_1\circ H$ is a
nonconstant Laurent polynomial,
and $G,G_1\in\C[x]$, so $H_1\circ H$ has no poles besides
$0$ and $\infty$.  But $H_1$ has poles at $0$ and $\infty$, so $H$
preserves $\{0,\infty\}$, and thus $H=\alpha x^s$ with $\alpha\in\C^*$
and $s\in\bZ$ (and $s\ne 0$).  If $s<0$ then, writing $\nu=1/x$, we have
$H_1\circ\nu=(x^n+1)/x^{n-r}$ and $\nu\circ H_2\circ\nu=H_2$ and
$G_2\circ\nu=(x+1)^n/x^{n-r}$, so by replacing $r$ by $n-r$ we again
have (\ref{Lunique}.2), but now with $H$ replaced by $x^{-s}/\alpha$.
Thus we may assume $s>0$, so the proof of Proposition~\ref{Lunique}
is complete.



\end{document}